\documentclass[11pt]{amsart}
\usepackage{geometry}                
\usepackage{graphicx}
\usepackage{amssymb}
\usepackage{epstopdf}
\usepackage{geometry} 
\usepackage{graphicx}
\usepackage{amsmath} 
\numberwithin{equation}{subsection}
\usepackage{commath}
\usepackage{enumitem}
\usepackage{changepage}
\usepackage{cleveref}
\usepackage[normalem]{ulem}
\usepackage{bbm}
\usepackage{subfig}
\usepackage{float}
\usepackage{url}
\usepackage{mathtools}
\usepackage{amsthm} 
\usepackage[disable]{todonotes} 
\DeclareMathOperator{\sgn}{sgn}
\DeclareMathOperator*{\esssup}{ess\,sup}
\theoremstyle{plain}
\newtheorem{theorem}{Theorem}[section]
\newtheorem*{theorem*}{Theorem}
\newtheorem{lemma}[theorem]{Lemma}

\newtheorem{proposition}[theorem]{Proposition}
\newenvironment{claim}[1]{\par\noindent\underline{Claim.}\space#1}{}
\newenvironment{claimproof}[1]{\par\noindent\emph{Proof of Claim.}\space#1}{\leavevmode\unskip\penalty9999 \hbox{}\nobreak\hfill\quad\hbox{$\blacksquare$}}
\theoremstyle{remark}
\newtheorem*{remark}{Remark}
\theoremstyle{definition}
\newtheorem{definition}[theorem]{Definition}
\allowdisplaybreaks
\title[Uniqueness with a Single Entropy Condition]{On Uniqueness of Solutions to Conservation Laws Verifying a Single Entropy Condition}
\author[Krupa]{Sam G. Krupa}
\address{Department of Mathematics\\ The University of Texas at Austin\\ Austin, TX 78712\\ USA}
\email{skrupa@math.utexas.edu}
\author[Vasseur]{Alexis F. Vasseur}
\address{Department of Mathematics\\ The University of Texas at Austin\\ Austin, TX 78712\\ USA}
\email{vasseur@math.utexas.edu}
\thanks{The second author was partly supported by NSF Grant DMS-1614918.}
\date{February 9th, 2018}                                           

\begin{document}
\begin{abstract}
For hyperbolic systems of conservation laws, uniqueness of solutions is still largely open. We aim to expand the theory of uniqueness for systems of conservation laws. One difficulty is that many systems have only one entropy. This contrasts with scalar conservation laws, where many entropies exist. It took until 1994 to show that one entropy is enough to ensure uniqueness of solutions for the scalar conservation laws (see Panov [{\em Mat. Zametki}, 55(5):116--129, 159, 1994]). This single entropy result was proven again by De Lellis, Otto and Westdickenberg about 10 years later [{\em Quart. Appl. Math.}, 62(4):687--700, 2004]. These two proofs both rely on the special connection between Hamilton--Jacobi equations and scalar conservation laws in one space dimension. However, this special connection does not extend to systems. In this paper, we prove the single entropy result for scalar conservation laws without using Hamilton--Jacobi.  Our proof lays out new techniques that are promising for showing uniqueness of solutions in the systems case.
\end{abstract}
\maketitle
\tableofcontents

\section{Introduction}
The present paper is concerned with the scalar conservation law in one space dimension:
\begin{align}
\label{conservation_law}
\begin{cases} 
   u_t+(A(u))_x=0\\
   u(x,0)=u^0(x).
  \end{cases}
\end{align}
In the scalar conservation law \eqref{conservation_law}, $u(x,t):\mathbb{R}\times[0,\infty)\to\mathbb{R}$ is the unknown, $u^0\in L^{\infty}(\mathbb{R})$ is the given initial data, and $A:\mathbb{R}\to\mathbb{R}$ is the given \emph{flux function}. In the present paper, we are concerned with $A\in C^2(\mathbb{R})$ strictly convex. A \emph{classical solution} of \eqref{conservation_law} is a locally Lipschitz function $u:\mathbb{R}\times[0,\infty)\to\mathbb{R}$ which satisfies $u_t+(A(u))_x=0$ almost everywhere and verifies $u(x,0)=u^0(x)$ for all $x\in\mathbb{R}$. We are also interested in \emph{weak solutions} of \eqref{conservation_law}. A weak solution to \eqref{conservation_law} is a locally bounded measurable function $u:\mathbb{R}\times[0,\infty)\to\mathbb{R}$ which satisfies 
\begin{align}
\int\limits_{0}^{\infty} \int\limits_{-\infty}^{\infty}[\partial_t \phi u +\partial_x \phi A(u)] \,dxdt + \int\limits_{-\infty}^{\infty} \phi(x,0)u^0(x)\,dx=0 \label{weak_solution}
\end{align}
for every Lipschitz continuous test function $\phi:\mathbb{R}\times[0,\infty)\to\mathbb{R}$, with compact support. In particular, every weak solution satisfies \eqref{conservation_law} in the sense of distributions. Note that a classical solution is also a weak solution.
\vskip0.3cm 
A pair of functions $\eta, q:\mathbb{R}\to\mathbb{R}$ are called an \emph{entropy} and \emph{entropy flux}, respectively, for the scalar conservation law \eqref{conservation_law} if
\begin{align}
q'(u)=\eta'(u)A'(u).
\end{align}
We say a weak solution $u$ of \eqref{conservation_law} is \emph{entropic} for the entropy $\eta$ if it satisfies the \emph{entropy inequality} 
\begin{align}
\label{entropy_inequality_distributional}
\partial_t \eta(u)+\partial_x q (u) \leq 0
\end{align}
in a distributional sense, where $q$ is any corresponding entropy flux. Precisely, 
\begin{align}
\int\limits_{0}^{\infty} \int\limits_{-\infty}^{\infty}[\partial_t \phi \eta(u) +\partial_x \phi q(u)] \,dxdt + \int\limits_{-\infty}^{\infty} \phi(x,0)\eta(u^0(x))\,dx\geq0 \label{entropy_inequality_integral_formulation}
\end{align}
for every nonnegative Lipschitz continuous test function $\phi:\mathbb{R}\times[0,\infty)\to\mathbb{R}$, with compact support. Kruzhkov \cite{MR0267257} proved existence and uniqueness for bounded weak solutions to \eqref{conservation_law} which are entropic for the large family of entropies $\{\eta_k\}_{k\in\mathbb{R}}$, where
\begin{align}
\eta_k(u)\coloneqq\abs{u-k}. \label{kruzkov_entropies}
\end{align}
For a bounded and measurable solution to \eqref{conservation_law}, being entropic for each of the $\eta_k$ is equivalent to being entropic for every convex entropy \cite[Proposition 2.3.4]{serre_book}. See Bolley, Brenier and Loeper \cite{brenier_optimal_transport76} for an extension of Kruzhkov's theory, based on the Wasserstein distance. Ole\u{\i}nik discovered \emph{``condition E,''} and proved existence and uniqueness for bounded weak solutions to \eqref{conservation_law} which satisfy it \cite{MR0094541}. A solution $u$ to \eqref{conservation_law} satisfies condition E if
\begin{align}
\label{condition_e}
\begin{cases}
\mbox{There exists a constant $C>0$ such that}\\
\hspace{1.3in}u(x+z,t)-u(x,t)\leq \frac{C}{t} z\\
\mbox{for all $t>0$, almost every $z>0$, and almost every $x\in\mathbb{R}$.}
\end{cases}
\end{align}
It is known that being entropic for each of the $\eta_k$ is equivalent to Ole\u{\i}nik's condition E when the conservation law \eqref{conservation_law} has a uniformly convex flux function $A$ \cite[p.~66 and p.~ 57]{serre_book}. When the solution to \eqref{conservation_law} is bounded and $A$ is strictly convex, we can assume $A$ is uniformly convex. 
\vskip0.3cm

For general systems, using the  $L^1$ theory, Bressan, Crasta, and Piccoli showed uniqueness in the class of solutions with small total variation in  \cite{MR1686652}. It can be interesting to study the uniqueness of the same solutions in a larger class. For example, for the $2\times 2$ Euler system, existence of solutions with large data is known. But uniqueness is still open for such solutions.

The present paper concerns the uniqueness of solutions to scalar conservation laws in one space dimension which are entropic for only one entropy.  However we are careful in our theory to develop techniques that we believe will extend to the systems case.

To attempt new progress on the theory of uniqueness, we take an entirely new approach. We use the method of relative entropy combined with the recent idea of stability up to a shift (first described by the second author in \cite{VASSEUR2008323}). Our methods are fundamentally $L^2$ theory.

The method of relative entropy allows for stability estimates to be made between a smooth solution to a conservation law and a weak solution entropic for at least one convex entropy. The method of relative entropy is powerful, applying to the cases of scalar, systems, and multiple space dimensions. 

The proof of stability estimates between weak and strong solutions relies on the fact that classical solutions of a hyperbolic system of conservation laws satisfy the entropy inequality \eqref{entropy_inequality_distributional} as an exact equality. Stability breaks down when a discontinuity is introduced into the classical solution. In particular, the relative entropy method is much more involved when considering shocks. For a first result, see DiPerna \cite{MR523630} for uniqueness. In the case of stability, the theory of stability up to a shift allows for discontinuities to be introduced into the smooth solution that the method of relative entropy considers. In this way, the method of relative entropy can be used to make comparisons between a weak solution and piecewise-Lipschitz solutions. This is how the uniqueness result of the present paper is proven.

Given the system \eqref{conservation_law} and entropy and entropy flux $\eta$ and $q$, respectively (or more generally, any hyperbolic system of conservation laws in multiple space dimensions endowed with any entropy), the method of relative entropy considers the quantity called the \emph{relative entropy}:
\begin{align}
\eta(a|b)\coloneqq \eta(a)-\eta(b)-\eta'(b)(a-b)
\end{align}
for all $a,b\in\mathbb{R}$. 

We have also the associated \emph{relative entropy-flux},
\begin{align}
q(a;b)\coloneqq q(a)-q(b)-\eta'(b)(A(a)-A(b))
\end{align}
for all $a,b\in\mathbb{R}$.

Both $\eta(a|b)$ and $q(a;b)$ are  locally quadratic in $a-b$. In particular, for all $a$ and $b$ in a fixed compact set, the strict convexity of $\eta\in C^2(\mathbb{R})$ gives
\begin{align}
\label{relative_entropy_controls_l2_4}
c^*(a-b)^2\leq \eta(a|b) \leq c^{**}(a-b)^2
\end{align}
for constants $c^*, c^{**}>0$. 

The method of relative entropy was invented by Dafermos \cite{doi:10.1080/01495737908962394,MR546634} and DiPerna \cite{MR523630} to prove weak-strong estimates: given a weak solution $u$ to \eqref{conservation_law}, entropic for $\eta$, and a classical solution $\bar{u}$, the method of relative entropy gives stability estimates on the growth in time of $\norm{u(\cdot,t)-\bar{u}(\cdot,t)}_{L^2}$ by considering the time derivative of $\int\eta(u|\bar{u})\,dx$ and using the entropy inequality \eqref{entropy_inequality_distributional} (in particular, see \Cref{extended_entropy_inequality_lemma}). Due to \eqref{relative_entropy_controls_l2_4}, the quantity $\eta(u|\bar{u})$ gives estimates of $L^2$-type, while being more amenable to study than the $L^2$ norm itself, due to the entropy inequality \eqref{entropy_inequality_distributional}. The relative entropy method is fundamentally an $L^2$ theory.

The aforementioned recent insight of the second author has allowed for discontinuities to exist in the otherwise smooth solution $\bar{u}$ that the method of relative entropy considers, while maintaining stability.
The key is that shocks have a contraction property in $L^2$ up to shift, even for large perturbations. For this, the discontinuities must not be allowed to move in time according to the conservation law, but instead their movement must be dictated by a special time-dependent function which shifts the solution. The first result in this program was by Leger \cite{Leger2011_original} for the scalar conservation law \eqref{conservation_law} for a strictly convex flux $A\in C^2(\mathbb{R})$. We now introduce the result of Leger. Let $u_L, u_R\in\mathbb{R}$ satisfy $u_L>u_R$. Let $\sigma$ satisfy $A(u_L)-A(u_R)=\sigma(u_L-u_R)$. Define
\begin{align}
  \phi (x)\coloneqq
  \begin{cases}
   u_L & \text{if } x<0 \\
   u_R & \text{if } x>0.
  \end{cases}
\end{align}
Let $u$ be any Kruzhkov solution. In this context, Leger proved the existence of a Lipschitz continuous function $h:[0,\infty)\to\mathbb{R}$ such that
\begin{align}
\norm{u(\cdot,t)-\phi(\cdot - h(t)-\sigma t)}_{L^2(\mathbb{R})}\leq\norm{u(\cdot,0)-\phi(\cdot)}_{L^2(\mathbb{R})}
\end{align} 
for all $t\geq0$. Note that by shifting the position of the shock wave as a function of time, all $L^2$ growth in time is killed. The shift function $h$ depends on $u$.  Leger gives control on $h$: $\abs{h(t)}\leq \lambda \norm{u(\cdot,0)-\phi(\cdot)}_{L^2(\mathbb{R})}\sqrt{t}$, for some constant $\lambda>0$. Leger only considered Kruzhkov solutions, but his methods are in fact very general and can be applied whenever a solution satisfies a strong trace property (see \Cref{strong_trace_condition} below) and is entropic for at least one strictly convex entropy $\eta\in C^2(\mathbb{R})$.

The second author and coworkers have been actively developing the theory of contraction up to a shift function. Progress has been made on systems of conservation laws in one space dimension by introducing the notion of a-contraction up to shift  \cite{MR3519973,MR3479527,MR3537479,serre_vasseur,Leger2011} and scalar viscous conservation laws in both one space dimension \cite{MR3592682} and multiple space dimensions \cite{multi_d_scalar_viscous_9122017}. For a more general overview of the theory of shifts and the relative entropy method in general, see \cite[Section 3-5]{MR3475284}. The theory of stability up to a shift has also been used to study asymptotic limits when the limit is discontinuous. See \cite{MR3333670} for the scalar case, and \cite{MR3421617} for the case of systems. There is a long history of using the  relative entropy method to study the asymptotic limit. However, results on the asymptotic limit which do not use shifts have only been able to consider the case when the limit function is Lipschitz continuous (see  \cite{MR1121850,MR1115587,MR1213991,MR2505730,MR2178222,MR2025302,MR1842343,MR1980855} and \cite{VASSEUR2008323} for a survey).

Current work only allows for a single discontinuity to exist in the otherwise smooth solution considered by the method of relative entropy. A natural progression is to try to allow for arbitrarily many discontinuities in the smooth solution, so that we can use the method of relative entropy to compare not just a smooth solution and a weak solution, but ideally two weak solutions. In the present paper, we use the method of relative entropy to compare a weak solution to a solution with arbitrarily many discontinuities. Thus, the present paper expands the theory of stability up to a shift.

The present paper lays out a new engine for proving uniqueness of solutions. In summary, the engine works by comparing two solutions using the method of relative entropy. One solution $u$ can be weak and entropic for just a single entropy. The second solution will have more regularity, but by allowing many discontinuities can still be from a dense class. By using the method of relative entropy to approximate the weak solution $u$ with a sequence of more regular solutions, we detect regularity in $u$. Hopefully, enough regularity in $u$ will ensure uniqueness.

The uniqueness result in this paper is for the scalar conservation laws. But, much of our work should generalize to systems. The method of relative entropy works for systems, and also multiple space dimensions. The theory of stability up to a shift has been developed for systems.  Many hyperbolic systems of conservation laws in one space dimension admit only one non-trivial entropy \cite[p.~238]{MR3475284}. In the present paper, we use only one entropy.

For the scalar conservation laws, uniqueness of solutions entropic for a single entropy is not new. The first proof was given by Panov in 1994 \cite{panov_uniquness}, for the system \eqref{conservation_law} with any flux $A\in C^2(\mathbb{R})$ strictly convex, any entropy $\eta\in C^1(\mathbb{R})$ strictly convex, and $u^0\in L^{\infty}(\mathbb{R})$. A second proof was given approximately 10 years later by De Lellis, Otto and Westdickenberg \cite{delellis_uniquneness}. Their result is stronger than Panov's result: their proof allows for various right-hand sides in the entropy inequality, and can consider unbounded functions \cite[p.~688]{delellis_uniquneness}.

However, both the proofs of Panov and De Lellis-Otto-Westdickenberg are fundamentally limited to the scalar conservation laws. Both proofs exploit the special connection between scalar conservation laws in one space dimension and Hamilton--Jacobi equations: the space derivative of the solution to a Hamilton--Jacobi equation is formally a solution to the associated scalar conservation law. It is well-known that the relation between scalar conservation laws in one space dimension and Hamilton--Jacobi equations breaks down in the more general case of hyperbolic systems of conservation laws in one space dimension (however, see \cite{jin_xin} for a formal connection between a general Hamilton--Jacobi equation in $n$ space dimensions and a \emph{weakly} hyperbolic system of conservation laws). 
\vskip0.3cm 
Before we can state precisely the uniqueness result proven in the present paper, let us introduce the strong trace property (originally introduced in \cite{Leger2011}):

\begin{definition}
\label{strong_trace_condition}
Let $u\in L^{\infty}(\mathbb{R}\times[0,\infty))$. Then $u$ verifies the \emph{strong trace property} if for any Lipschitz continuous function $h: [0,\infty)\to\mathbb{R}$, there exist $u_{+},u_{-}\in L^{\infty}([0,\infty))$ such that
\begin{align}
\lim_{n\to\infty} \int\limits_{0}^{T}\esssup_{y\in(0,\frac{1}{n})}\abs{u(h(t)+y,t)-u_{+}(t)}\,dt=\lim_{n\to\infty} \int\limits_{0}^{T}\esssup_{y\in(-\frac{1}{n},0)}\abs{u(h(t)+y,t)-u_{-}(t)}\,dt=0
\end{align}
for all $T>0$.
\end{definition}

Note that any function $u\in L^{\infty}(\mathbb{R}\times[0,\infty))$ will satisfy the strong trace property if $u$ has a representative such that for any fixed $h$, the right and left limits  
\begin{align}
\lim_{y\to0^{+}} u(h(t)+y,t) \hspace{.5in}\text{and}\hspace{.5in} \lim_{y\to0^{-}} u(h(t)+y,t) 
\end{align}
exist for almost every t. In particular, any function $u\in L^{\infty}(\mathbb{R}\times[0,\infty))$ with a representative in $BV_{\text{loc}}$ will satisfy the strong trace property. But the strong trace property is weaker than $BV_{\text{loc}}$. 

The result we show in this paper is

\begin{theorem}[Main theorem]

\label{main_theorem}
 Let $u\in L^{\infty}(\mathbb{R}\times[0,\infty))$ be a weak solution  with initial data $u^0\in L^{\infty}(\mathbb{R})$ to the scalar conservation law in one space dimension \eqref{conservation_law} with flux $A\in C^2(\mathbb{R})$ strictly convex. Assume $u$ satisfies the entropy inequality \eqref{entropy_inequality_distributional} for at least one strictly convex entropy $\eta\in C^2(\mathbb{R})$. Further, assume $u$ satisfies the strong trace property (\Cref{strong_trace_condition}).
 
 Then $u$ is the unique solution to \eqref{conservation_law} verifying \eqref{condition_e} and with initial data $u^0$.
\end{theorem}

We briefly outline the present paper and the proof of the main theorem, \Cref{main_theorem}. Let $u$ be any bounded weak solution to \eqref{conservation_law} with initial data $u^0\in L^{\infty}(\mathbb{R})$, entropic for a strictly convex entropy $\eta\in C^2(\mathbb{R})$ and satisfying the strong trace property (\Cref{strong_trace_condition}). In the scalar case, the special structure we are trying to detect is \eqref{condition_e}. To prove $u$ verifies \eqref{condition_e}, we use shift functions very similar to the ones constructed by Leger \cite{Leger2011_original} and in \cite{serre_vasseur} to give, at each fixed time $T$, a sequence of piecewise Lipschitz continuous functions $\{\psi_\epsilon\}_{\epsilon>0}$ defined on a subset of the real line which converge in $L^2$ to $u(\cdot,T)$ as $\epsilon\to0^{+}$. The $\psi_\epsilon$ will be constructed by gluing together at time $T$ various classical solutions to \eqref{conservation_law}, each of which satisfy \eqref{condition_e}. Thus, the $\psi_\epsilon$ will be shown to verify \eqref{condition_e}. By classical measure theory, the $\psi_\epsilon$ converge (up to a subsequence) to $u(\cdot,T)$ pointwise almost everywhere. Thus, $u$ will have the same structure as the $\psi_\epsilon$ and will satisfy \eqref{condition_e} for time $T$, where $T$ is arbitrary. This will complete the proof of \Cref{main_theorem}.

In \Cref{Preliminaries}, we give some preliminary results and facts which we will need. Then, in \Cref{construction_of_shift} we give a construction of a shift function as suited to our needs in this paper, based on the construction in  \cite{serre_vasseur}. With the shift construction out of the way, in \Cref{main_prop_section} we state and prove the main proposition (\Cref{main_lemma}) which gives the existence of the $\psi_\epsilon$. The proof of the main proposition is broken up into two lemmas (\Cref{psi_exists_with_limits} and \Cref{LIWAS_dense}). To conclude, in \Cref{main_prop_implies_main_theorem} we let $\epsilon\to0^{+}$ (up to a subsequence) to explain how the main proposition implies the main theorem (\Cref{main_theorem}).

\section{Preliminaries}\label{Preliminaries}

Throughout this paper, we will assume that there is a constant $B>0$ such that all of the weak solutions $w$ to \eqref{conservation_law} that we consider satisfy 
\begin{align}
\label{define_B}
\norm{w}_{L^{\infty}(\mathbb{R}\times[0,\infty))}\leq B.
\end{align}

For a function $u\in L^{\infty}(\mathbb{R}\times[0,\infty))$ that verifies the strong trace property (\Cref{strong_trace_condition}), and a Lipschitz continuous function $h:[0,\infty)\to\mathbb{R}$, we use the notation $u(h(t)\pm,t)$ to denote the values $u_{\pm}(t)$ at a time $t$ when $u$ has a left and right trace according to \Cref{strong_trace_condition}, where $u_{\pm}$ are as in the definition of strong trace. For time $t$ when $u$ has a strong trace, the values $u_{\pm}(t)$ are well-defined, and hence $u(h(t)\pm,t)$ are also well-defined.

Most of the things in this preliminary section are well-known. But we include these facts here and their proofs so as to make sure we are not using more than we think. In \Cref{main_theorem} we do not assume that the solution $u$ to the scalar conservation law is Kruzhkov.

\begin{lemma}[Properties of solutions with Lipschitz initial data]\label{lip_initial_data_properties}
Let $v_i^0\in L^{\infty}(\mathbb{R})$ be nondecreasing Lipschitz continuous functions for $i=1,2$, verifying $v_1^0(x)\geq v_2^0(x)$ for all $x\in\mathbb{R}$. Let $v_i$ denote the unique solution to \eqref{conservation_law} satisfying \eqref{condition_e} and with initial data $v_i^0$. 
Then the following holds:
\begin{enumerate}[label=(\alph*)]
\item \label{lip_initial_1}
The $v_i$ are classical solutions to \eqref{conservation_law} on $\mathbb{R}\times[0,\infty)$.

\item \label{lip_initial_2}
The $v_i$ are given by the method of characteristics. In other words, for each $(x,t)\in\mathbb{R}\times[0,\infty)$, there exists $x^0\in\mathbb{R}$ such that $v_1(x,t)=v_1^0(x^0)$ and $x=x^0+tA'(v_1^0(x^0))$ (and similarly for $v_2$). 

\item \label{lip_initial_6}
The $v_i$ satisfy \eqref{condition_e} with $C=1/ \inf A''$.

\item \label{lip_initial_5}
For $i=1,2$ and $t\geq0$, $v_i(\cdot,t)$ is a nondecreasing function in the $x$ variable.

\item \label{lip_initial_3}
$\norm{v_i^0}_{L^{\infty}(\mathbb{R})}=\norm{v_i}_{L^{\infty}(\mathbb{R}\times[0,\infty))}$ for $i=1,2$.

\item \label{lip_initial_4}
$v_1(x,t)\geq v_2(x,t)$ for all $(x,t)\in\mathbb{R}\times[0,\infty)$.
\end{enumerate}
\begin{proof}
From \cite[p.~176-7]{dafermos_big_book}, we know that for each $i$ there exists a classical solution $v_i$ to \eqref{conservation_law} on $\mathbb{R}\times[0,\infty)$ with initial data $v_i^0$ which is given by the method of characteristics. 

For a fixed $i$, we now check by direct method that $v_i$ satisfies \eqref{condition_e}. By uniqueness of solutions satisfying \eqref{condition_e}, this will prove parts \ref{lip_initial_1} and \ref{lip_initial_2}.

Fix $t,z>0$ and $x\in\mathbb{R}$. Because $v_i$ is given by the method of characteristics, there exists an $x^0$ and $\tilde{x}^0$ such that
\begin{align}
\begin{cases}
v_i(x,t)=v_i^0(x^0) \\
x=x^0+tA'(v_i^0(x^0))
\end{cases}
\end{align}
and
\begin{align}
\begin{cases}
v_i(x+z,t)=v_i^0(\tilde{x}^0) \\
x+z=\tilde{x}^0+tA'(v_i^0(\tilde{x}^0)).
\end{cases}
\end{align}

We have
\begin{equation}
\begin{aligned}\label{characteristics1}
v_i(x,t)=(A')^{-1}\bigg(\frac{x-{x}^0}{t}\bigg)\\
v_i(x+z,t)=(A')^{-1}\bigg(\frac{x+z-\tilde{x}^0}{t}\bigg)\\
\end{aligned}
\end{equation}
where the functional inverse $(A')^{-1}$ exists because $A$ is strictly convex.

Because the characteristics do not cross, we have
\begin{align}\label{characteristics4}
0 < \tilde{x}^0 -x^0=z+t(A'(v_i^0(x^0))-A'(v_i^0(\tilde{x}^0)))<z.
\end{align}

Then from \eqref{characteristics4},
\begin{align}\label{characteristics3}
0<\frac{x+z-\tilde{x}^0}{t}-\frac{x-{x}^0}{t}<\frac{z}{t}.
\end{align}
Finally,
\begin{equation}
\begin{aligned}\label{characteristics2}
&\inf A'' \bigg[(A')^{-1}\bigg(\frac{x+z-\tilde{x}^0}{t}\bigg)-(A')^{-1}\bigg(\frac{x-{x}^0}{t}\bigg)\bigg] \\
&\leq A'((A')^{-1}\bigg(\frac{x-{x}^0}{t}\bigg))-A'((A')^{-1}\bigg(\frac{x+z-\tilde{x}^0}{t}\bigg))<\frac{z}{t}.
\end{aligned}
\end{equation}

Lines \eqref{characteristics1}, \eqref{characteristics3} and \eqref{characteristics2} imply $v_i$ satisfies \eqref{condition_e}. In particular, note that we can take $C=1/ \inf A''$ in \eqref{condition_e}.

This proves parts \ref{lip_initial_1}, \ref{lip_initial_2}, and \ref{lip_initial_6}.

Part \ref{lip_initial_5} follows immediately from part \ref{lip_initial_2}. Similarly, part \ref{lip_initial_3} follows immediately from parts \ref{lip_initial_1} and \ref{lip_initial_2}.

We now show how part \ref{lip_initial_4} follows from part \ref{lip_initial_2}. We argue by contradiction. Assume there exists $(x,t)\in\mathbb{R}\times[0,\infty)$ such that $v_1(x,t)<v_2(x,t)$.

Then we have just proven that $v_1$ and $v_2$ are given by the method of characteristics. Thus there exists $x^0$ and $\tilde{x}^0$ such that
\begin{equation}
\begin{aligned}\label{characteristics_a1}
x&=x^0+tA'(v_1^0(x^0))\\
&=\tilde{x}^0+tA'(v_2^0(\tilde{x}^0))
\end{aligned}
\end{equation}
and
\begin{align}
\label{characteristics_a2}
v_1(x,t)=v_1^0(x^0)<v_2^0(\tilde{x}^0)=v_2(x,t).
\end{align}

Then from \eqref{characteristics_a1},
\begin{align}
\label{characteristics_a3}
x^0-\tilde{x}^0=t(A'(v_2^0(\tilde{x}^0))-A'(v_1^0(x^0))).
\end{align}

Then the right-hand side of \eqref{characteristics_a3} is nonnegative, which means $x^0\geq \tilde{x}^0$. Thus,
\begin{align}
v_1^0(x^0)\geq v_1^0(\tilde{x}^0)\geq v_2^0(\tilde{x}^0)
\end{align}
because $v_1^0(x)\geq v_2^0(x)$ for all $x\in\mathbb{R}$. However, this gives a contradiction with \eqref{characteristics_a2}. This proves part \ref{lip_initial_4}.
\end{proof}
\end{lemma}

\begin{lemma}[An extension of the entropy inequality]
\label{extended_entropy_inequality_lemma}
Let $\bar{u}^0\in L^{\infty}(\mathbb{R})$ be a Lipschitz continuous and nondecreasing function. Let $\bar{u}$ be the unique solution to \eqref{conservation_law} with initial data $\bar{u}^0$ and which satisfies \eqref{condition_e}. Let $u\in L^{\infty}(\mathbb{R}\times[0,\infty))$ be a weak solution to \eqref{conservation_law} with initial data $u^0$. Assume that $u$ is entropic for the strictly convex entropy $\eta\in C^2(\mathbb{R})$. 

Then
\begin{align}
\label{relative_entropy_inequality_distributional}
\partial_t \eta (u|\bar{u})+\partial_x q(u;\bar{u}) \leq 0
\end{align}
in the sense of distributions. 

In other words, the following holds for all nonnegative, Lipschitz continuous test functions $\phi$ on $\mathbb{R}\times [0,\infty)$ with compact support:
\begin{align}
\label{extended_entropy_inequality}
\int\limits_{0}^{\infty} \int\limits_{-\infty}^{\infty}[\partial_t \phi \eta(u|\bar{u})+\partial_x \phi q(u;\bar{u})] \,dxdt + \int\limits_{-\infty}^{\infty} \phi(x,0)\eta(u^0(x)|\bar{u}^0(x))\,dx\geq0.
\end{align}
\begin{remark}
The inequality \eqref{relative_entropy_inequality_distributional} extends the entropy inequality \eqref{entropy_inequality_distributional}.
\end{remark}
\begin{proof}
By part \ref{lip_initial_1} of \Cref{lip_initial_data_properties}, $\bar{u}$ is a classical solution. Then, \Cref{extended_entropy_inequality_lemma} follows immediately from the following inequality:
\begin{align}
\label{fact_from_dafermos}
\begin{split}
\int\limits_{0}^{\infty} \int\limits_{-\infty}^{\infty}[\partial_t \phi \eta(u|\bar{u})+\partial_x \phi q(u;\bar{u})] \,dxdt + \int\limits_{-\infty}^{\infty} \phi(x,0)\eta(u^0(x)|\bar{u}^0(x))\,dx\geq \\
\int\limits_{0}^{\infty} \int\limits_{-\infty}^{\infty}\phi\partial_x \bar{u} \eta''(\bar{u})[A(u)-A(\bar{u})-A'(\bar{u})(u-\bar{u})]\,dxdt. 
\end{split}
\end{align}

Dafermos gives this inequality as one of the central steps in the proof of  weak-strong stability. See \cite[p.~124, line (5.2.10)]{dafermos_big_book}.

By part \ref{lip_initial_5} of \Cref{lip_initial_data_properties}, $\bar{u}$ is increasing in $x$. Furthermore, we have $\eta'',A''>0$. Thus the right-hand side of \eqref{fact_from_dafermos} is controlled from below by zero.
\end{proof}
\end{lemma}

\begin{lemma}\label{rigorous_dissipation_rate_calc_middle_part}
Let $[t^*,t^{**})$ be a bounded interval. Let $h_1(t), h_2(t): [t^*,t^{**})\to \mathbb{R}$ be Lipschitz continuous functions, such that $h_2(t)-h_1(t)>0$ for all $t\in[t^*,t^{**})$.

Let $\bar{u}^0\in L^{\infty}(\mathbb{R})$ be a Lipschitz continuous and nondecreasing function. Let $\bar{u}$ be the unique solution to \eqref{conservation_law} with initial data $\bar{u}^0$ and which satisfies \eqref{condition_e}. Let $u\in L^{\infty}(\mathbb{R}\times[0,\infty))$ be a weak solution to \eqref{conservation_law} with initial data $u^0$. Assume $u$ is entropic for the strictly convex entropy $\eta\in C^2(\mathbb{R})$. Furthermore,  assume that $u$ verifies the strong trace property (\Cref{strong_trace_condition}).

Then for almost every $a,b\in[t^*,t^{**})$ with $a<b$,
\begin{equation}
\begin{aligned}\label{rigorous_dissipation_rate_calc_middle_part_ineq}
&\int\limits_{h_1(b)}^{h_2(b)}\eta (u(x,b)|\bar{u}(x,b))\,dx-\int\limits_{h_1(a)}^{h_2(a)}\eta (u(x,a)|\bar{u}(x,a))\,dx\leq\\
&\hspace{.2in}\int\limits_{a}^{b}\bigg[ q(u(h_1(t)+,t);\bar{u}(h_1(t),t))-q(u(h_2(t)-,t);\bar{u}(h_2(t),t))+\\&\hspace{.2in}\dot{h}_2(t)\eta(u(h_2(t)-,t)|\bar{u}(h_2(t),t))-\dot{h}_1(t)\eta(u(h_1(t)+,t)|\bar{u}(h_1(t),t))\bigg]\,dt.
\end{aligned}
\end{equation}

If $t^*=0$, then \eqref{rigorous_dissipation_rate_calc_middle_part_ineq} holds for $a=0$ and almost every $b\in(t^*,t^{**})$.

\begin{proof} 
For $0<\epsilon<\min\{b-a,t^{**}-b\}$, define 

\begin{align}
 \chi_\epsilon (t)\coloneqq
  \begin{cases}
   0 & \text{if } t<a \\
   \frac{1}{\epsilon}(t-a) & \text{if } a\leq t<a+\epsilon \\
   1 & \text{if } a+\epsilon\leq t\leq b\\
   -\frac{1}{\epsilon}(t-(b+\epsilon)) & \text{if } b<t \leq b+\epsilon\\
   0 & \text{if } b+\epsilon<t.
  \end{cases}
\end{align}

Let $\delta\coloneqq \inf_{t\in[a,b+\epsilon]}h_2(t)-h_1(t)$. Note $\delta>0$. Then for $0<\epsilon<\frac{\delta}{2}$, define

\begin{align}
 \psi_{\epsilon}(x,t)\coloneqq
  \begin{cases}
   0 & \text{if } x<h_1(t) \\
   \frac{1}{\epsilon}(x-h_1(t)) & \text{if } h_1(t)<x<h_1(t)+\epsilon \\
   1 & \text{if } h_1(t)+\epsilon<x<h_2(t)-\epsilon\\
   -\frac{1}{\epsilon}(x-h_2(t)) & \text{if } h_2(t)-\epsilon<x<h_2(t)\\
   0 & \text{if } h_2(t)<x.
  \end{cases}
\end{align}

The $\psi_{\epsilon}(x,t)$ and $\chi_{\epsilon} (t)$ are from \cite[p.~765]{Leger2011_original}.

Use $\psi_{\epsilon}(x,t) \chi_{\epsilon}(t)$ as a test function for \eqref{extended_entropy_inequality}.
The result is
\begin{equation}
\begin{aligned}
\label{before_rearranged_dissipation_rate}
&\frac{1}{\epsilon}\int\limits_{a}^{b}\bigg[\int\limits_{h_1(t)}^{h_1(t)+\epsilon} q(u;\bar{u})\,dx-\int\limits_{h_2(t)-\epsilon}^{h_2(t)} q(u;\bar{u})\,dx+\int\limits_{h_2(t)-\epsilon}^{h_2(t)} \dot{h}_2(t)\eta(u|\bar{u})\,dx-\int\limits_{h_1(t)}^{h_1(t)+\epsilon}\dot{h}_1(t)\eta(u|\bar{u})\,dx\bigg]\,dt \\
&\hspace{.2in}+ \frac{1}{\epsilon}\int\limits_{a}^{a+\epsilon} \int\limits_{h_1(t)}^{h_2(t)}\eta(u|\bar{u}) \,dxdt -\frac{1}{\epsilon}\int\limits_{b}^{b+\epsilon} \int\limits_{h_1(t)}^{h_2(t)}\eta(u|\bar{u}) \,dxdt+O(\epsilon)\geq0.
\end{aligned}
\end{equation}

Let $\epsilon\to0^{+}$ in \eqref{before_rearranged_dissipation_rate}. Use the dominated convergence theorem, the Lebesgue differentiation theorem, and the fact that $u$ verifies the strong trace property (\Cref{strong_trace_condition}). This gives \eqref{rigorous_dissipation_rate_calc_middle_part_ineq} for almost every $a,b\in[t^*,t^{**})$ with $a<b$.

When $t^*=0$, we want to show that  \eqref{rigorous_dissipation_rate_calc_middle_part_ineq} holds for $a=t^*=0$. This follows because we include the boundary term corresponding to $t=0$ in \eqref{entropy_inequality_integral_formulation}, and hence the boundary term corresponding to $t=0$ appears also in \eqref{extended_entropy_inequality}.

We now show that when $t^*=0$, \eqref{rigorous_dissipation_rate_calc_middle_part_ineq} holds for $a=t^*=0$.

For $a=t^*=0$ and $0<\epsilon<t^{**}-b$, define 

\begin{align}
 \chi_\epsilon^0 (t)\coloneqq
  \begin{cases}
   1 & \text{if } 0\leq t< b\\
   -\frac{1}{\epsilon}(t-(b+\epsilon)) & \text{if } b\leq t < b+\epsilon\\
   0 & \text{if } b+\epsilon \leq t.
  \end{cases}
\end{align}

The $\chi_\epsilon^0$ is borrowed from \cite[p.~124]{dafermos_big_book}.  Test \eqref{extended_entropy_inequality} with $\psi_{\epsilon}(x,t) \chi_{\epsilon}^0(t)$. This gives

\begin{equation}
\begin{aligned}
\label{before_rearranged_dissipation_rate_zero}
&\frac{1}{\epsilon}\int\limits_{a}^{b}\bigg[\int\limits_{h_1(t)}^{h_1(t)+\epsilon} q(u;\bar{u})\,dx-\int\limits_{h_2(t)-\epsilon}^{h_2(t)} q(u;\bar{u})\,dx+\int\limits_{h_2(t)-\epsilon}^{h_2(t)} \dot{h}_2(t)\eta(u|\bar{u})\,dx-\int\limits_{h_1(t)}^{h_1(t)+\epsilon}\dot{h}_1(t)\eta(u|\bar{u})\,dx\bigg]\,dt \\
&\hspace{1in}+ \int\limits_{h_1(0)}^{h_2(0)}\eta(u^0(x)|\bar{u}^0(x)) \,dx -\frac{1}{\epsilon}\int\limits_{b}^{b+\epsilon} \int\limits_{h_1(t)}^{h_2(t)}\eta(u|\bar{u}) \,dxdt+O(\epsilon)\geq0.
\end{aligned}
\end{equation}

Finally, let $\epsilon\to0^{+}$ in \eqref{before_rearranged_dissipation_rate_zero}. Once again, invoke the dominated convergence theorem, the Lebesgue differentiation theorem, and use that $u$ verifies the strong trace property (\Cref{strong_trace_condition}). We receive \eqref{rigorous_dissipation_rate_calc_middle_part_ineq} for $a=0$ and for almost every $b\in(t^*,t^{**})$.
\end{proof}
\end{lemma}

\begin{lemma}\label{left_right_ap_limits}
Let $\bar{u}^0\in L^{\infty}(\mathbb{R})$ be a Lipschitz continuous and nondecreasing function. Let $\bar{u}$ be the unique solution to \eqref{conservation_law} with initial data $\bar{u}^0$ and which satisfies \eqref{condition_e}. Let $u\in L^{\infty}(\mathbb{R}\times[0,\infty))$ be a weak solution to \eqref{conservation_law} with initial data $u^0$. Assume that $u$ is entropic for at least one strictly convex entropy $\eta\in C^2(\mathbb{R})$. Assume also that $u$ verifies the strong trace property (\Cref{strong_trace_condition}). 

Then for all $c,d\in\mathbb{R}$ verifying $c<d$, the approximate right- and left-hand limits
\begin{align}\label{ap_right_left_limits_exist_44}
&{\rm ap}\,\lim_{t\to {t_0}^{\pm}}\int\limits_{c}^{d}\eta (u(x,t)|\bar{u}(x,t))\,dx
\end{align}
exist for all $t_0\in(0,\infty)$ and verify
\begin{align}
\label{left_and_right_limits_order}
&{\rm ap}\,\lim_{t\to {t_0}^{-}}\int\limits_{c}^{d}\eta (u(x,t)|\bar{u}(x,t))\,dx\geq {\rm ap}\,\lim_{t\to {t_0}^{+}}\int\limits_{c}^{d}\eta (u(x,t)|\bar{u}(x,t))\,dx.
\end{align}
The approximate right-hand limit also exists for $t_0=0$ and verifies 
\begin{align}
\label{right_limit_at_zero}
\int\limits_{c}^{d}\eta (u^0(x)|\bar{u}^0(x))\,dx\geq {\rm ap}\,\lim_{t\to {0}^{+}}\int\limits_{c}^{d}\eta (u(x,t)|\bar{u}(x,t))\,dx.
\end{align}
\begin{proof}
For some constant $C>0$ to be chosen momentarily, define the function $\Gamma:[0,\infty)\to\mathbb{R}$, 
\begin{align}
\Gamma(t)\coloneqq \int\limits_{c}^{d}\eta (u(x,t)|\bar{u}(x,t))\,dx -Ct.
\end{align}
Apply \eqref{rigorous_dissipation_rate_calc_middle_part_ineq} to the case when $h_1(t)=c$ and $h_2(t)=d$ for all $t$. The integrand on the right-hand side of \eqref{rigorous_dissipation_rate_calc_middle_part_ineq} is bounded. Thus, there exists some constant $C>0$ such that $\Gamma(t)\geq\Gamma(s)$ for almost every $t$ and $s$ verifying $t<s$. This means that there exists a function which agrees with $\Gamma$ almost everywhere and is non-increasing. This implies that $\Gamma$ has approximate left and right  limits. In particular, we conclude that the approximate right- and left-hand limits \eqref{ap_right_left_limits_exist_44}
exist for all $t_0\in(0,\infty)$ and verify \eqref{left_and_right_limits_order}.
Note the approximate right-hand limit also exists for $t_0=0$ and because \eqref{rigorous_dissipation_rate_calc_middle_part_ineq} holds for the time $a=0$, the approximate right-hand limit verifies \eqref{right_limit_at_zero} at time zero.

\end{proof}
\end{lemma}

\section{Construction of the shift} \label{construction_of_shift}
In this section, we present a proof of
\begin{proposition}[Existence of the shift function]
\label{shift_theorem}
Let $u\in L^{\infty}(\mathbb{R}\times[0,\infty))$ be a weak solution to \eqref{conservation_law}, entropic for at least one strictly convex entropy $\eta\in C^2(\mathbb{R})$. Assume also that $u$ verifies the strong trace property (\Cref{strong_trace_condition}). Let $\bar{u}^0_i\in L^{\infty}(\mathbb{R})$ be Lipschitz continuous and nondecreasing functions for $i=1,2$, verifying $\bar{u}^0_1(x)\geq\bar{u}_2^0(x)$ for all $x\in\mathbb{R}$. Let $\bar{u}_i\in L^{\infty}(\mathbb{R}\times[0,\infty))$ be solutions  to \eqref{conservation_law} verifying \eqref{condition_e} and with initial data $\bar{u}^0_i$, for $i=1,2$. Let $t^*\geq0$ be some fixed time, and let $x_0\in\mathbb{R}$ be some fixed space coordinate.

Then for $\epsilon>0$ there exists a Lipschitz continuous function $h_\epsilon:[t^*,\infty)\to\mathbb{R}$ such that
\begin{align}
h_\epsilon(t^*)=x_0, \hspace{.3in} \mbox{\upshape Lip}[h_\epsilon]\leq  \sup_{\hspace{.02in}\abs{x}\leq\norm{u}_{L^{\infty}}} \abs{A'}
\end{align}
and
\begin{align}
q(u_{+};\bar{u}_2)-q(u_{-};\bar{u}_1)-\dot{h}_\epsilon(\eta(u_{+}|\bar{u}_2)-\eta(u_{-}|\bar{u}_1))\leq \epsilon
\end{align}
for almost every $t\in[t^*,\infty)$, where $u_{\pm}=u(h_\epsilon(t)\pm,t)$ and $\bar{u}_i=\bar{u}_i(h_\epsilon(t),t)$. 
\end{proposition}

In \cite{serre_vasseur}, the authors build a general framework for the construction of the shift functions necessary for $L^2$ stability in the general case of systems. They also apply their framework to the scalar case, recovering the result of Leger \cite{Leger2011_original}.

We present here a construction of the shift function based on the work \cite{serre_vasseur}. We modify the construction slightly. In \cite{serre_vasseur}, the data $\bar{u}_i$ are constant. In contrast, we assume the $\bar{u}_i$ are Lipschitz continuous solutions to \eqref{conservation_law}. The presentation is also simplified due to our focus only on the scalar case. Lastly, in the general framework for systems in \cite{serre_vasseur}, there is a  ``constraint'' labeled with equation number 7 \cite[p.~4]{serre_vasseur}. Although this constraint is used in \cite{serre_vasseur} for the scalar case in particular, we find that the proof for scalar in \cite{serre_vasseur} as-is does not require it.

Lastly, the paper \cite{serre_vasseur} considers only Kruzhkov's solutions \cite[p.~9]{serre_vasseur} which are entropic for the family of entropies $\{\abs{u-k}\}_{k\in\mathbb{R}}$. We find that the proofs in \cite{serre_vasseur} hold without modification for solutions which are not necessarily of this type.

\subsection{Lemmas necessary for the proof of \Cref{shift_theorem}}

We need the following structural lemmas.

\begin{lemma}[Structural lemma on entropic shocks]
\label{entropic_shocks}
Consider the equation \eqref{conservation_law}, endowed with a strictly convex entropy $\eta\in C^2(\mathbb{R})$, and an associated entropy flux $q$. Let $u_L,u_R,\sigma\in\mathbb{R}$ satisfy
\begin{align}
A(u_R)-A(u_L)=\sigma(u_R-u_L).
\end{align}
Then 
\begin{align}
\label{entropic_shock}
q(u_R)-q(u_L)\leq \sigma (\eta(u_R)-\eta(u_L))
\end{align}
if and only if $u_L\geq u_R$. Which is to say, the shock $(u_L,u_R,\sigma)$ is entropic for the entropy $\eta$ if and only if $u_L\geq u_R$.
\begin{remark}
For the scalar conservation law \eqref{conservation_law} with flux $A$ strictly convex, we know that the ``physical'' or ``admissible'' shock wave solutions will be the ones with left- and right-hand states $u_L$ and $u_R$, respectively, which satisfy $u_L>u_R$. This is the \emph{Lax entropy condition} for the scalar case. It is immediate that shock solutions satisfying \eqref{condition_e} verify the Lax entropy condition. \Cref{entropic_shocks} says that shock solutions to the scalar conservation law \eqref{conservation_law} endowed with the entropy $\eta$, also verify the Lax entropy condition.
\end{remark}

\begin{proof}
We only consider the case when $u_L\neq u_R$.

Let
\begin{align}
\Lambda \coloneqq q(u_R)-q(u_L)-\sigma (\eta(u_R)-\eta(u_L))= \int\limits_{u_L}^{u_R} A'(u)\eta'(u)-\sigma\eta'(u)\,du.
\end{align}
We want to show $\Lambda \leq 0$ if and only if $u_L>u_R$.

We can write 
\begin{align}
\sigma=\frac{A(u_R)-A(u_L)}{u_R-u_L}=\frac{1}{u_R-u_L}\int\limits_{u_L}^{u_R} A'(v)\,dv.
\end{align}
Thus 
\begin{align}
\Lambda
&= 
\int\limits_{u_L}^{u_R} \eta'(u)\bigg(A'(u)- \frac{1}{u_R-u_L}\int\limits_{u_L}^{u_R} A'(v)\,dv\bigg)\,du\\
&=
 \int\limits_{u_L}^{u_R} \eta'(u)\bigg(\frac{1}{u_R-u_L}\int\limits_{u_L}^{u_R} A'(u)-A'(v)\,dv\bigg)\,du\\
 &=
 \frac{1}{u_R-u_L}\int\limits_{u_L}^{u_R}\int\limits_{u_L}^{u_R}\eta'(u)(A'(u)-A'(v))\,dvdu\\
 &=
 \frac{1}{u_R-u_L}\bigg(\frac{1}{2}\iint\limits_{D}\eta'(u)(A'(u)-A'(v))\,dvdu-\frac{1}{2}\iint\limits_{D}\eta'(v)(A'(u)-A'(v))\,dvdu\bigg)
\end{align}
where $D\coloneqq I(u_L,u_R)\times I(u_L,u_R)$. We use $I(a,b)$ to denote the interval with endpoints $a$ and $b$.

Finally,
\begin{align}
\Lambda=\frac{1}{2}\frac{1}{u_R-u_L}\iint\limits_{D} (\eta'(u)-\eta'(v))(A'(u)-A'(v))\,dvdu.
\end{align}
Then, by strict convexity of $A$ and $\eta$, the quantity 
\begin{align}
\label{always_nonnegative}
\iint\limits_{D} (\eta'(u)-\eta'(v))(A'(u)-A'(v))\,dvdu
\end{align}
is always positive. 
Thus, the sign of $\Lambda$ is given by the sign of $u_R-u_L$. This completes the proof.
\end{proof}
\end{lemma}

\begin{lemma}[Structural lemma from \cite{serre_vasseur}]
\label{d_sm_rh_lemma}
Let $u_L,u_R,u_{-},u_{+}\in\mathbb{R}$ satisfy $u_L > u_R$ and $u_{-} > u_{+}$. Let $\eta\in C^2(\mathbb{R})$ be a strictly convex entropy, with associated entropy flux $q$. Define
\begin{align}
\sigma(u_{-},u_{+})\coloneqq \frac{A(u_{+})-A(u_{-})}{u_{+}-u_{-}}.\label{rh_velocity}
\end{align}

Then, 
\begin{align}
q(u_{+};u_R)-q(u_{-};u_L)-\sigma(u_{-},u_{+})(\eta(u_{+}|u_R)-\eta(u_{-}|u_L)) \leq 0. \label{d_rh}
\end{align}

Moreover, if
\begin{align}
\eta(u|u_L)=\eta(u|u_R) \label{denom_zero}
\end{align}
for some $u\in\mathbb{R}$, then
\begin{align}
q(u;u_R)-q(u;u_L)< 0. \label{d_sm}
\end{align}
\end{lemma}

\Cref{d_sm_rh_lemma} was proven in \cite{serre_vasseur}. For completeness, we give the proof of this lemma in the appendix (\Cref{d_sm_rh_lemma_proof}). 

\begin{lemma}[Structural lemma]
\label{v_epsilon_bounded}
Let $\eta\in C^2(\mathbb{R})$ be a strictly convex entropy for \eqref{conservation_law}. Let $q$ be the associated entropy flux. For $\epsilon>0$, define the function $V_\epsilon: \{(u,u_L,u_R)\in \mathbb{R}^3 | u_L\geq u_R\}\to \mathbb{R}$:
\begin{align}
  V_\epsilon(u,u_L,u_R)\coloneqq
  \begin{cases}
   \frac{[q(u;u_R)-q(u;u_L)-\epsilon]_+}{\eta(u|u_R)-\eta(u|u_L)} & \text{if } \eta(u|u_R)\neq\eta(u|u_L) \\
   0       & \text{if } \eta(u|u_R)=\eta(u|u_L),
  \end{cases}
\end{align}
where $[\hspace{.025in}\cdot\hspace{.025in}]_{+}\coloneqq \max(0,\cdot)$.

Then, $V_\epsilon$ verifies
\begin{align}
\abs{V_\epsilon(u,u_L,u_R)} \leq \abs{A'(u)}. \label{final_part_v_epsilon_bounded}
\end{align} 
\begin{proof}
Let 
\[
  V(u,u_L,u_R)\coloneqq
  \begin{cases}
   \frac{[q(u;u_R)-q(u;u_L)]_+}{\eta(u|u_R)-\eta(u|u_L)} & \text{if } \eta(u|u_R)\neq\eta(u|u_L) \\
   0       & \text{if } \eta(u|u_R)=\eta(u|u_L).
  \end{cases}
\]

In order to show \eqref{final_part_v_epsilon_bounded}, the idea is to show
\begin{align}
\abs{V(u,u_L,u_R)} \leq \abs{A'(u)}\label{v_control}
\end{align}
and then use the basic inequality 
\[\abs{V_\epsilon(u,u_L,u_R)}\leq \abs{V(u,u_L,u_R)}.\]

The proof of \eqref{v_control} depends on controlling the quantity 
\begin{align}\label{quantity_we_want_to_control}
 \frac{q(u;u_R)-q(u;u_L)}{\eta(u|u_R)-\eta(u|u_L)}
 \end{align}
 for the $(u,u_L,u_R)$ values that make $q(u;u_R)-q(u;u_L)\geq0$. 
 
We have three cases where we get control on the quantity \eqref{quantity_we_want_to_control}: $u>u_L>u_R$,  $u_L>u_R>u$, and $u_R\leq u\leq u_L$. We begin with some elementary facts which will be used repeatedly. 

Remark that 
\begin{align}
\partial_b \eta(a|b)=\eta''(b)(b-a) \mbox{\hspace{.3in}and\hspace{.3in}} \partial_b q(a;b)=\eta''(b)(A(b)-A(a)).
\end{align}

In particular, we can write for all $a,b,c\in\mathbb{R}$,
\begin{align}
q(a;b)-q(a;c)= \int\limits_{c}^{b}\eta''(v)(A(v)-A(a))\,dv= \int\limits_{c}^{b}\eta''(v)(v-a)\frac{A(v)-A(a)}{v-a}\,dv.
\end{align}

We can now begin the casework to prove \Cref{v_epsilon_bounded}.

\emph{Case} $u>u_L>u_R$

We note that for all $u,u_L,u_R\in\mathbb{R}$ such that $u>u_L>u_R$,
\begin{align}
\eta(u|u_R)-\eta(u|u_L)=\int\limits_{u_L}^{u_R}\partial_v \eta(u|v)\,dv=\int\limits_{u_L}^{u_R}\eta''(v)(v-u)\,dv>0
\end{align}
by strict convexity of $\eta$.
We conclude
\begin{align}
\label{sign_denom_1}
\eta(u|u_R)-\eta(u|u_L)>0.
\end{align}

For all $u,u_L,u_R\in\mathbb{R}$ such that $u>u_L>u_R$ and $\eta(u|u_R)-\eta(u|u_L)\neq0$, we compute 
\begin{align}
\frac{q(u;u_R)-q(u;u_L)}{\eta(u|u_R)-\eta(u|u_L)}
&=
\frac{\int\limits_{u_L}^{u_R}\eta''(v)(v-u)\frac{A(v)-A(u)}{v-u}\,dv}{\eta(u|u_R)-\eta(u|u_L)} \label{refer_back_1}
\shortintertext{By the mean value theorem and strict convexity of $A$,  $\frac{A(v)-A(u)}{v-u}<A'(u)$ for $v\in[u_R,u_L]$. Thus, recalling also the strict convexity of $\eta$ and \eqref{sign_denom_1}, we continue from \eqref{refer_back_1} to get}
&<
\frac{A'(u)\int\limits_{u_L}^{u_R}\eta''(v)(v-u)\,dv}{\eta(u|u_R)-\eta(u|u_L)}
=A'(u).
\end{align}

In summary,

\begin{align}\label{bound_v_case1}
\frac{q(u;u_R)-q(u;u_L)}{\eta(u|u_R)-\eta(u|u_L)}<A'(u).
\end{align}

\emph{Case} $u_L>u_R>u$

Our method is analogous to the case $u>u_L>u_R$ above.

Remark that for all $u,u_L,u_R\in\mathbb{R}$ such that $u_L>u_R>u$,
\begin{align}
\eta(u|u_R)-\eta(u|u_L)=\int\limits_{u_L}^{u_R}\partial_v \eta(u|v)\,dv=\int\limits_{u_L}^{u_R}\eta''(v)(v-u)\,dv<0
\end{align}
by strict convexity of $\eta$.
We conclude
\begin{align}
\label{sign_denom_2}
\eta(u|u_R)-\eta(u|u_L)<0.
\end{align}

For all $u,u_L,u_R\in\mathbb{R}$ such that $u_L>u_R>u$ and $\eta(u|u_R)-\eta(u|u_L)\neq0$, we calculate 
\begin{align}
\frac{q(u;u_R)-q(u;u_L)}{\eta(u|u_R)-\eta(u|u_L)}
&=
\frac{\int\limits_{u_L}^{u_R}\eta''(v)(v-u)\frac{A(v)-A(u)}{v-u}\,dv}{\eta(u|u_R)-\eta(u|u_L)}\label{refer_back_2}
\shortintertext{By the mean value theorem and strict convexity of $A$,  $\frac{A(v)-A(u)}{v-u}>A'(u)$ for $v\in[u_R,u_L]$. Thus, recalling also the strict convexity of $\eta$ and \eqref{sign_denom_2}, we continue from \eqref{refer_back_2} to get}
&>
\frac{A'(u)\int\limits_{u_L}^{u_R}\eta''(v)(v-u)\,dv}{\eta(u|u_R)-\eta(u|u_L)}
=A'(u).
\end{align}

To summarize,

\begin{align}\label{bound_v_case2}
\frac{q(u;u_R)-q(u;u_L)}{\eta(u|u_R)-\eta(u|u_L)}>A'(u).
\end{align}

\emph{Case} $u_R\leq u\leq u_L$

This case is slightly different than the two cases above.

For all $u,u_L,u_R\in\mathbb{R}$ such that $u_R\leq u\leq u_L$ and $\eta(u|u_R)-\eta(u|u_L)\neq0$,
\begin{align}
 \frac{q(u;u_R)-q(u;u_L)}{\abs{\eta(u|u_R)-\eta(u|u_L)}}
 &=
 \frac{\int\limits_{u_L}^{u_R}\eta''(v)(v-u)\frac{A(v)-A(u)}{v-u}\,dv}{\abs{\eta(u|u_R)-\eta(u|u_L)}}\\
&=
 \frac{\int\limits_{u_L}^{u}\eta''(v)(v-u)\frac{A(v)-A(u)}{v-u}\,dv+\int\limits_{u}^{u_R}\eta''(v)(v-u)\frac{A(v)-A(u)}{v-u}\,dv}{\abs{\eta(u|u_R)-\eta(u|u_L)}}
 \shortintertext{By the mean value theorem and strict convexity of $A$, $\frac{A(v)-A(u)}{v-u}>A'(u)$ for $v\in(u,u_L]$, and $\frac{A(v)-A(u)}{v-u}<A'(u)$ for $v\in[u_R,u)$. Thus, recalling also the strict convexity of $\eta$,}
 &<
 \frac{A'(u)\int\limits_{u_L}^{u}\eta''(v)(v-u)\,dv+A'(u)\int\limits_{u}^{u_R}\eta''(v)(v-u)\,dv}{\abs{\eta(u|u_R)-\eta(u|u_L)}}\\
 &=
A'(u)\sgn(\eta(u|u_R)-\eta(u|u_L)).
\end{align}

We receive,
\begin{align}\label{bound_v_case3}
 \frac{q(u;u_R)-q(u;u_L)}{\abs{\eta(u|u_R)-\eta(u|u_L)}}<A'(u)\sgn(\eta(u|u_R)-\eta(u|u_L)).
\end{align}

We combine all the cases $u>u_L>u_R$,  $u_L>u_R>u$, and $u_R\leq u\leq u_L$, and in particular \eqref{sign_denom_1}, \eqref{bound_v_case1}, \eqref{sign_denom_2}, \eqref{bound_v_case2}, and \eqref{bound_v_case3}. We keep in mind that we only consider $(u,u_L,u_R)$ values that make $q(u;u_R)-q(u;u_L)\geq0$.

In conclusion,
\[\abs{V_\epsilon(u,u_L,u_R)}\leq \abs{V(u,u_L,u_R)} \leq \abs{A'(u)}.\]
\end{proof}
\end{lemma}

\subsection{Proof of \Cref{shift_theorem}}

This proof is based on the work \cite[p.~7-8]{serre_vasseur}. We modify the proof to consider $\bar{u}_i$ which are non-constant. 

Define $V_\epsilon: \{(u,u_L,u_R)\in \mathbb{R}^3 | u_L\geq u_R\}\to \mathbb{R}$:
\begin{align}
  V_\epsilon(u,u_L,u_R)\coloneqq
  \begin{cases}
   \frac{[q(u;u_R)-q(u;u_L)-\epsilon]_+}{\eta(u|u_R)-\eta(u|u_L)} & \text{if } \eta(u|u_R)\neq\eta(u|u_L) \\
   0       & \text{if } \eta(u|u_R)=\eta(u|u_L),
  \end{cases}
\end{align}
where $[\hspace{.025in}\cdot\hspace{.025in}]_{+}\coloneqq \max(0,\cdot)$.

The function $V_\epsilon$ is continuous. For $(u,u_L,u_R)$ such that $\eta(u|u_R)\neq\eta(u|u_L)$, $V_\epsilon$ is clearly continuous. By \eqref{d_sm}, $V_\epsilon=0$ on some neighborhood of $\{(u,u_L,u_R)\in \mathbb{R}^3 | \eta(u|u_R)=\eta(u|u_L)\}$. Thus, $V_\epsilon$ is continuous.

Remark that by part \ref{lip_initial_1} of \Cref{lip_initial_data_properties}, the $\bar{u}_i$ are classical solutions, and by part \ref{lip_initial_4} of \Cref{lip_initial_data_properties}, $\bar{u}_1(x,t)\geq\bar{u}_2(x,t)$ for all $(x,t)\in\mathbb{R}\times[0,\infty)$.

We construct a solution to 
\begin{align}
\label{filippov}
\begin{cases} 
  \dot{h}_\epsilon=V_\epsilon(u(h_\epsilon(t),t),\bar{u}_1(h_\epsilon(t),t),\bar{u}_2(h_\epsilon(t),t))\\
   h_\epsilon(t^*)=x_0
  \end{cases}
\end{align}
in the Filippov sense. We use the following lemma:

\begin{lemma}
\label{Filippov20}
There exists a Lipschitz function $h_\epsilon:[t^*,\infty)\to\mathbb{R}$ such that
\begin{align}
&\hspace{-1.5in}h_\epsilon(t^*)=x_0,\label{Filippov1}\\
&\hspace{-1.5in}\norm{\dot{h}_\epsilon}_{L^{\infty}}\leq\norm{V_\epsilon}_{L^{\infty}},\label{Filippov2}\\
\shortintertext{and}
\dot{h}_\epsilon(t)\in I[V_\epsilon(u_{-},\bar{u}_1,\bar{u}_2),V_\epsilon(u_{+},\bar{u}_1,\bar{u}_2)],\label{Filippov3}
\end{align}
for almost every $t>0$, where $u_{\pm}=u(h_\epsilon(t)\pm,t)$, and $\bar{u}_i=\bar{u}_i(h_\epsilon(t),t)$ for $i=1,2$. Here, $I[a,b]$ denotes the closed interval with endpoints $a$ and $b$.

Moreover, for almost every $t>0$,
\begin{align}
A(u_{+})-A(u_{-})=\dot{h}_\epsilon(u_{+}-u_{-}),\label{Filippov4}\\
q(u_{+})-q(u_{-})\leq\dot{h}_\epsilon (\eta(u_{+})-\eta(u_{-})),\label{Filippov5}
\end{align}
which means that for almost every $t>0$, either the shock $(u_{-},u_{+},\dot{h}_\epsilon)$ is an entropic discontinuity for the entropy $\eta$ or $u_{-}=u_{+}$.
\end{lemma}

The proof of \eqref{Filippov1}, \eqref{Filippov2} and \eqref{Filippov3} is nearly identical to the proof of Proposition 1 in \cite{Leger2011}. For completeness, we provide a proof of \eqref{Filippov1}, \eqref{Filippov2} and \eqref{Filippov3} in the appendix (\Cref{proof_of_Filippov}). The properties \eqref{Filippov4} and \eqref{Filippov5} in fact hold for any Lipschitz function $h: [0,\infty)\to\mathbb{R}$. These properties are well-known in the $BV$ case. When $u$ only satisfies the strong trace property (\Cref{strong_trace_condition}), \eqref{Filippov4} and \eqref{Filippov5} are given in Lemma 6 in \cite{Leger2011}. We do not include a proof of \eqref{Filippov4} and \eqref{Filippov5} here; a proof is given in the appendix in \cite{Leger2011}.

We are now ready to prove \Cref{shift_theorem}.
\begin{proof}[Proof of \Cref{shift_theorem}]
For almost every $t$ such that $u_{-}=u_{+}$, we borrow notation from \cite{serre_vasseur} and write  $u_{\pm}\coloneqq u_{-}=u_{+}$, and then by \Cref{Filippov20} we have $\dot{h}_\epsilon=V_\epsilon(u_{\pm},\bar{u}_1,\bar{u}_2)$. This gives, 
\begin{align}
q(u_{\pm};\bar{u}_2)-q(u_{\pm};\bar{u}_1)-\dot{h}_\epsilon(\eta(u_{\pm}|\bar{u}_2)-\eta(u_{\pm}|\bar{u}_1)) \leq \epsilon.
\end{align}

For almost every $t$ such that $u_{-}\neq u_{+}$, \Cref{Filippov20} says $\dot{h}_\epsilon=\sigma(u_{-},u_{+})$, where $\sigma(u_{-},u_{+})=(A(u_{+})-A(u_{-}))/(u_{+}-u_{-})$. Then,

\begin{align}
&q(u_{+};\bar{u}_2)-q(u_{-};\bar{u}_1)-\dot{h}_\epsilon(\eta(u_{+}|\bar{u}_2)-\eta(u_{-}|\bar{u}_1)) =\\
&q(u_{+};u_R)-q(u_{-};u_L)-\sigma(u_{-},u_{+})(\eta(u_{+}|u_R)-\eta(u_{-}|u_L)), \\
\shortintertext{and then by \eqref{d_rh}, \Cref{entropic_shocks}, and \Cref{Filippov20} again,}
&\leq 0.
\end{align}

Thus, for almost every $t\in[t^*,\infty)$
\begin{align}
q(u_{+};\bar{u}_2)-q(u_{-};\bar{u}_1)-\dot{h}_\epsilon(\eta(u_{+}|\bar{u}_2)-\eta(u_{-}|\bar{u}_1))\leq \epsilon.
\end{align}

Finally, by \Cref{v_epsilon_bounded} and \Cref{Filippov20},
\begin{align}
\mbox{Lip}[h_\epsilon]\leq  \sup_{\hspace{.02in}\abs{x}\leq\norm{u}_{L^{\infty}}} \abs{A'}.
\end{align}
\end{proof}

\section{Main proposition} \label{main_prop_section}
The proof of \Cref{main_theorem} will follow from 
\begin{proposition}[Main proposition]

\label{main_lemma}
Let $u\in L^{\infty}(\mathbb{R}\times[0,\infty))$ be a weak solution to \eqref{conservation_law}, with initial data $u^0$. Let $\eta\in C^2(\mathbb{R})$ be a strictly convex entropy. Assume that $u$ is entropic for the entropy $\eta$ and verifies the strong trace property (\Cref{strong_trace_condition}).

Then for all $R, T, \epsilon>0$, there exists a function $\psi:[-R,R] \to\mathbb{R}$ verifying:
\begin{enumerate}[label=(\alph*)]
\item \label{relative_entropy_stable}
\begin{align*}
\int\limits_{\abs{x}\leq R}\eta(u(x,T)|\psi(x))\,dx \leq \epsilon.
\end{align*}
\item \label{psi_satisfies_condition_e}
\begin{align*}
\psi(x+z)-\psi(x)\leq \frac{c}{T}z
\end{align*}
for $x\in[-R,R]$ and $z>0$ with $x+z\in[-R,R]$ and where $c=1/\inf{A''}$.
\item \label{psi_bounded_u}
\begin{align*}
\norm{\psi}_{L^\infty([-R,R])}\leq \norm{u}_{L^\infty(\mathbb{R}\times[0,\infty))}.
\end{align*}
\end{enumerate}
\end{proposition}


We decompose the proof of \Cref{main_lemma} into two lemmas. We utilize functions of the form
\begin{align}
\hspace{1in}
\label{LIWAS}
\begin{cases}
\intertext{For $v\in L^{\infty}(\mathbb{R})$, there exists a finite set of $x_i$ with}
\hspace{1in}-\infty=x_0<x_1<x_2<\cdots<x_N<x_{N+1}=\infty
\intertext{such that $v$ is nondecreasing and Lipschitz continuous on $(x_i,x_{i+1})$}
\intertext{for $i=0,\ldots,N$, and}
\hspace{1.7in}\lim\limits_{\substack{x\to x_i \\ x<x_i}} v(x) \geq \lim\limits_{\substack{x\to x_i \\ x>x_i}} v(x) \\
\mbox{for $1\leq i \leq N$.}
\end{cases}
\end{align}

\begin{lemma}
\label{psi_exists_with_limits}
Let $R,T>0$. Let $u\in L^{\infty}(\mathbb{R}\times[0,\infty))$ be a weak solution to \eqref{conservation_law} with initial data $u^0$. Assume $u$ is entropic for at least one strictly convex entropy $\eta\in C^2(\mathbb{R})$. Assume that $u$ verifies the strong trace property (\Cref{strong_trace_condition}). Choose $s>0$ to verify $\abs{q(a;b)}\leq s \eta(a|b)$ for all $a,b\in[-B,B]$, where $B$ is defined in \eqref{define_B}.

Then for all $N\in\mathbb{N}\cup\{0\}$, we have:

For any $v_1^0,\ldots,v_{N+1}^0\in L^{\infty}(\mathbb{R})$  Lipschitz continuous nondecreasing functions verifying $v_i^0(x)\geq v_{i+1}^0(x)$ for all $1\leq i \leq N$ and $x\in\mathbb{R}$, and for any $t^*\in[0,T]$ and real numbers $x_0,\ldots,x_{N+1}$ verifying $x_0=-R+(t^*-T)s<x_1<\cdots<x_{N}<R-(t^*-T)s=x_{N+1}$ the following holds:

There exist $N+2$ real numbers $x_{0,T},\ldots,x_{N+1,T}$ such that $x_{0,T}=-R\leq x_{1,T}\leq \cdots\leq x_{N,T}\leq R=x_{N+1,T}$ and
 
\begin{align}
\label{l2_stability1}
{\rm ap}\,\lim_{t\to {T}^{+}}\sum_{i=0}^{N} \int\limits_{x_{i,T}}^{x_{i+1,T}}\eta(u(x,t)|v_{i+1}(x,t))\,dx
\leq 
{\rm ap}\,\lim_{t\to {t^*}^{+}}\sum_{i=0}^{N} \int\limits_{x_i}^{x_{i+1}}\eta(u(x,t)|v_{i+1}(x,t))\,dx
\end{align}
where $v_i$ is the unique solution to \eqref{conservation_law} with initial data $v_i^0$ and verifying \eqref{condition_e}.
\end{lemma}

\begin{lemma}[Density of functions of form \eqref{LIWAS} in $L^2$]
\label{LIWAS_dense}
Let $M,\epsilon>0$. Then for all $f\in L^2([-M,M])$, there is a function  $f_\epsilon: \mathbb{R}\to\mathbb{R}$ of the form \eqref{LIWAS} such that 
\begin{align}
\norm{f-f_\epsilon}_{L^2([-M,M])} < \epsilon, \\
\norm{f_\epsilon}_{L^\infty(\mathbb{R})} \leq \norm{f}_{L^\infty([-M,M])},\label{l_infinity_control_on_LIWAS_dense}
\end{align}
and all of the discontinuities of $f_\epsilon$ are contained in $(-M,M)$.
\end{lemma}

\subsection{\Cref{psi_exists_with_limits} and \Cref{LIWAS_dense} imply main proposition}

As in \Cref{psi_exists_with_limits}, we choose $s>0$ such that $\abs{q(a;b)}\leq s \eta(a|b)$ for all $a,b\in[-B,B]$, where $B$ is defined in \eqref{define_B}. We also choose $c^{**}>0$ such that 
\begin{align}
\label{control_rel_entropy_by_above}
\eta(a|b)\leq c^{**}(a-b)^2
\end{align}
for all $a,b\in[-B,B]$.

By \Cref{LIWAS_dense}, there exists a function $v^0\in L^{\infty}(\mathbb{R})$ of the form \eqref{LIWAS} such that
\begin{align}
\label{l2_distance_LIWAS_guy}
\norm{u^0-v^0}_{L^2([-R-sT,R+sT])} < \sqrt{\frac{\epsilon}{c^{**}}}
\end{align}
and if there is at least one discontinuity in $v^0$, the discontinuities are at points $x_1<x_2<\cdots<x_N$ for some $N\in\mathbb{N}$, and where $x_i\in(-R-sT,R+sT)$ for  all $1\leq i \leq N$.
 
If $v^0$ contains at least one discontinuity, define the functions $v_i^0 : \mathbb{R}\to\mathbb{R}$ for $1\leq i \leq N+1$ as follows:

\[
  v_1^0 (x)\coloneqq
  \begin{cases}
   v^0(x) & \text{if } x<x_1 \\
   \operatorname*{sup}_{x_1< y < x}\max(v^0(x_1-),v^0(y)) & \text{if } x_1<x 
  \end{cases}
\]

For $2\leq i \leq N$,
\[
  v_i^0 (x)\coloneqq
  \begin{cases}
  \operatorname*{inf}_{x < y < x_{i-1}}\min(v^0(x_{i-1}+),v_{i-1}^0(y)) & \text{if } x<x_{i-1} \\
   v^0(x) & \text{if } x_{i-1}<x<x_{i} \\
   \operatorname*{sup}_{x_{i}< y< x}\max(v^0(x_{i}-),v^0(y)) & \text{if } x_{i}<x
  \end{cases}
\]

And

\[
  v_{N+1}^0 (x)\coloneqq
  \begin{cases}
  \operatorname*{inf}_{x < y < x_{N}}\min(v^0(x_{N}+),v_{N}^0(y)) & \text{if } x<x_{N} \\
   v^0(x) & \text{if } x_{N}<x.
  \end{cases}
\]

If $v^0$ has no discontinuities, then $N=0$ and we define $v^0_1\coloneqq v^0$. 

By construction, the $v^0_i$ are Lipschitz continuous, nondecreasing in $x$, and verify $v^0_i(x)\geq v^0_{i+1}(x)$ for all $x\in\mathbb{R}$ and $1\leq i \leq N$. We also have 
\begin{align}\label{sup_norm_control_on_v0i}
\norm{v^0_i}_{L^\infty(\mathbb{R})}\leq\norm{v^0}_{L^\infty(\mathbb{R})}
\end{align}
for $1\leq i \leq N+1$.

Let $v_i$ denote the unique solution to \eqref{conservation_law} with initial data $v^0_i$ and which satisfies \eqref{condition_e}.

From \Cref{psi_exists_with_limits}, we have $N+2$ real numbers $x_{0,T},\ldots,x_{N+1,T}$ such that $x_{0,T}=-R\leq x_{1,T}\leq \cdots\leq x_{N,T}\leq R=x_{N+1,T}$ and
 
\begin{align}
\label{l2_stability_proof_of_main_prop}
{\rm ap}\,\lim_{t\to {T}^{+}}\sum_{i=0}^{N} \int\limits_{x_{i,T}}^{x_{i+1,T}}\eta(u(x,t)|v_{i+1}(x,t))\,dx
\leq 
{\rm ap}\,\lim_{t\to {0}^{+}}\sum_{i=0}^{N} \int\limits_{x_i}^{x_{i+1}}\eta(u(x,t)|v_{i+1}(x,t))\,dx
\end{align}
where $x_0\coloneqq-R-sT$ and $x_{N+1}\coloneqq R+sT$.

We now control the right-hand side of \eqref{l2_stability_proof_of_main_prop}. Recalling \Cref{left_right_ap_limits}, we have
\begin{equation}
\begin{aligned}\label{control_rhs_l2_stability_proof_of_main_prop}
{\rm ap}\,\lim_{t\to {0}^{+}}\sum_{i=0}^{N} \int\limits_{x_i}^{x_{i+1}}\eta(u(x,t)|v_{i+1}(x,t))\,dx
&\leq 
\sum_{i=0}^{N} \int\limits_{x_i}^{x_{i+1}}\eta(u^0(x)|v_{i+1}^0(x))\,dx
\shortintertext{Then, from the definition of the $v_i^0$,}
&=\int\limits_{-R-sT}^{R+sT}\eta(u^0(x)|v^0(x))\,dx
\shortintertext{Using \eqref{control_rel_entropy_by_above}, \eqref{l_infinity_control_on_LIWAS_dense}, and that $\norm{u^0}_{L^\infty(\mathbb{R})}\leq \norm{u}_{L^{\infty}(\mathbb{R}\times[0,\infty))}$,}
&\leq c^{**} \int\limits_{-R-sT}^{R+sT}(u^0(x)-v^0(x))^2\,dx
\shortintertext{Then, from \eqref{l2_distance_LIWAS_guy}}
&< \epsilon.
\end{aligned}
\end{equation}

On the other hand, by the convexity of $\eta$,
\begin{align}\label{control_lhs_l2_stability_proof_of_main_prop}
\sum_{i=0}^{N} \int\limits_{x_{i,T}}^{x_{i+1,T}}\eta(u(x,T)|v_{i+1}(x,T))\,dx
\leq {\rm ap}\,\lim_{t\to {T}^{+}}\sum_{i=0}^{N} \int\limits_{x_{i,T}}^{x_{i+1,T}}\eta(u(x,t)|v_{i+1}(x,t))\,dx.
\end{align}

Combining \eqref{l2_stability_proof_of_main_prop}, \eqref{control_rhs_l2_stability_proof_of_main_prop} and \eqref{control_lhs_l2_stability_proof_of_main_prop}, we find
\begin{align}\label{combine_rhs_lhs_l2_stability_proof_of_main_prop}
\sum_{i=0}^{N} \int\limits_{x_{i,T}}^{x_{i+1,T}}\eta(u(x,T)|v_{i+1}(x,T))\,dx < \epsilon.
\end{align}

We define $\psi:[-R,R] \to\mathbb{R}$,

\begin{align}
\psi(x)\coloneqq
   \begin{cases}
   v_1(x,T) & \text{if } -R < x< x_{1,T} \\
   v_2(x,T) & \text{if } x_{1,T} < x< x_{2,T} \\
   \hspace{.24in} \vdots \\
   v_{N+1}(x,T) & \text{if } x_{N,T} < x< R.
   \end{cases}
\end{align}

By \eqref{combine_rhs_lhs_l2_stability_proof_of_main_prop}, $\psi$ satisfies part \ref{relative_entropy_stable} of \Cref{main_lemma}. 

We now show $\psi$ satisfies part \ref{psi_satisfies_condition_e} of \Cref{main_lemma}: this follows because each of the $v_i$ satisfy \eqref{condition_e}. In particular,
\begin{align}\label{condition_e_for_v_i}
v_i(x+z,T)-v_i(x,T)\leq \frac{C}{T} z
\end{align}
for all $z>0$ and all $x\in\mathbb{R}$. From part \ref{lip_initial_6} of \Cref{lip_initial_data_properties}, we can take $C=1/\inf A''$ in \eqref{condition_e_for_v_i}. Further, by part \ref{lip_initial_4} of \Cref{lip_initial_data_properties}, $v_i(x,T)\geq v_{i+1}(x,T)$ for all $x\in\mathbb{R}$ and $1\leq i \leq N$. This gives part \ref{psi_satisfies_condition_e} of \Cref{main_lemma}.

Part \ref{psi_bounded_u} of \Cref{main_lemma} follows from \eqref{l_infinity_control_on_LIWAS_dense}, \eqref{sup_norm_control_on_v0i}, and part \ref{lip_initial_3} of \Cref{lip_initial_data_properties}.

This completes the proof of \Cref{main_lemma}.

\subsection{Proof of \Cref{psi_exists_with_limits}}

We prove this lemma by strong induction on $N$.

\uline{Base case}

For $N=0$, let $v_1^0\in L^{\infty}(\mathbb{R})$ be any Lipschitz continuous nondecreasing function. Let $v_1$ be the unique solution to \eqref{conservation_law} with initial data $v_1^0$ and verifying \eqref{condition_e}.

Let $h_0(t)\coloneqq-R+(t-T)s$ and $h_{1}(t)\coloneqq R-(t-T)s$. Then, from \Cref{rigorous_dissipation_rate_calc_middle_part}, \Cref{left_right_ap_limits}, and the dominated convergence theorem,
\begin{equation}
\begin{aligned}\label{dissipation_calculation_base}
&{\rm ap}\,\lim_{t\to {T}^{+}} \int\limits_{h_{0}(T)}^{h_{1}(T)}\eta(u(x,t)|v_1(x,t))\,dx- {\rm ap}\,\lim_{t\to {t^{*}}^{+}}\int\limits_{h_{0}(t^{*})}^{h_{1}(t^{*})}\eta(u(x,t)|v_{1}(x,t))\,dx\\
&\leq
\int\limits_{t^{*}}^{T}\bigg[ q(u(h_0(t)+,t);v_1(h_0(t),t))-q(u(h_{1}(t)-,t);v_1(h_{1}(t),t))\\
&\hspace{.2in}+\dot{h}_{1}(t)\eta(u(h_{1}(t)-,t)|v_1(h_{1}(t),t))-\dot{h}_0(t)\eta(u(h_0(t)+,t)|v_1(h_0(t),t))\bigg]\,dt\\
&=\int\limits_{t^{*}}^{T}q(u(h_0(t)+,t);v_1(h_0(t),t))-s\eta(u(h_0(t)+,t)|v_1(h_0(t),t))\,dt\\
&+\int\limits_{t^{*}}^{T}-q(u(h_1(t)-,t);v_1(h_1(t),t))-s\eta(u(h_1(t)-,t)|v_1(h_1(t),t))\,dt\\
&\leq 0
\end{aligned}
\end{equation}
by the definition of $h_0,h_2$ and $s$.

We get,
\begin{align}
{\rm ap}\,\lim_{t\to {T}^{+}} \int\limits_{-R}^{R}\eta(u(x,t)|v_1(x,t))\,dx \leq {\rm ap}\,\lim_{t\to {t^{*}}^{+}}\int\limits_{-R+(t^*-T)s}^{R-(t^*-T)s}\eta(u(x,t)|v_{1}(x,t))\,dx.
\end{align}

Thus \Cref{psi_exists_with_limits} holds for the base case $N=0$.

\uline{Induction step}

Suppose that $K\in\mathbb{N}\cup\{0\}$ is given such that \Cref{psi_exists_with_limits} holds for $N=0,1,\ldots, K$.

We now prove that \Cref{psi_exists_with_limits} holds for $N=K+1$. Let $v_1^0,v_2,^0\ldots,v_{K+2}^0\in L^{\infty}(\mathbb{R})$ be any Lipschitz continuous nondecreasing functions, satisfying $v_i^0(x)\geq v_{i+1}^0(x)$ for all $1\leq i \leq K+1$ and $x\in\mathbb{R}$. For $1\leq i \leq K+2$, let  $v_i$ be the unique solution to \eqref{conservation_law} with initial data $v_i^0$ and verifying \eqref{condition_e}. Let $t^*\in[0,T]$ and $-R+(t^*-T)s<x_1<\cdots<x_{K+1}<R-(t^*-T)s$ be arbitrary.

Let $\epsilon>0$. By \Cref{shift_theorem} we can construct Lipschitz continuous functions $h_{\epsilon,1},\ldots,\\h_{\epsilon,K+1}$ on the interval $[t^*,T]$ such that for $1\leq i \leq K+1$, $h_{\epsilon,i}(t^*)=x_i$ and
\begin{align}
\label{dissipation_rate3}
q(u^{i}_{+};v_{i+1})-q(u^{i}_{-};v_i)-\dot{h}_{\epsilon,i}(\eta(u^{i}_{+}|v_{i+1})-\eta(u^{i}_{-}|v_i))\leq \frac{\epsilon}{T(K+1)}
\end{align}
for almost every $t\in[t^*,T]$, where $u^{i}_{\pm}=u(h_{\epsilon,i}(t)\pm,t)$ and $v_l=v_l(h_{\epsilon,i}(t),t)$ for $l=i,i+1$. 

To simplify the exposition, denote 
\begin{align}
h_{\epsilon,0}(t)\coloneqq -R+(t-T)s,\\
h_{\epsilon,K+2}(t)\coloneqq R-(t-T)s.
\end{align}

The Lipschitz constants of the $h_{\epsilon,i}$ are uniformly bounded in $\epsilon$ by \Cref{shift_theorem}. Thus, by Arzel\`a--Ascoli there exists a sequence $\{\epsilon_j\}_{j\in\mathbb{N}}$ such that $\epsilon_j\to0^{+}$ and for each $0\leq i \leq K+2$, $h_{\epsilon_j,i}$ converges uniformly on $[t^*,T]$ to a Lipschitz function $h_{i}$.

Then, let $t^{**}$ be the first time that there exist two of the $h_{i}$ for $0\leq i \leq K+2$ that are equal to each other. If such a time does not exist, let $t^{**}\coloneqq T$.  

We now show:

\begin{claim}
\begin{align}
\label{relative_entropy_decreasing_claim1}
{\rm ap}\,\lim_{t\to {t^{**}}^{+}}\sum_{i=0}^{K+1} \int\limits_{h_{i}(t^{**})}^{h_{i+1}(t^{**})}\eta(u(x,t)|v_{i+1}(x,t))\,dx
&\leq 
{\rm ap}\,\lim_{t\to {t^{*}}^{+}}\sum_{i=0}^{K+1} \int\limits_{h_{i}(t^{*})}^{h_{i+1}(t^{*})}\eta(u(x,t)|v_{i+1}(x,t))\,dx. 
\end{align}
\end{claim}
\begin{claimproof}

Due to the uniform convergence of the $h_{\epsilon_j,i}$ as $j\to\infty$, for each $\tau\in[t^*,t^{**})$, there exists $J_\tau>0$ large enough such that $h_{\epsilon_j,i+1}(t)- h_{\epsilon_j,i}(t)>0$ for all $t\in[t^*,\tau]$, $j>J_\tau$ and $0\leq i \leq K+1$.

Then, for almost every $t$ and $\tau$ verifying $t^*\leq t < \tau < t^{**}$, we have from  \Cref{rigorous_dissipation_rate_calc_middle_part} (for $j>J_\tau$):
\begin{align*}
&\sum_{i=0}^{K+1} \int\limits_{h_{\epsilon_j,i}(\tau)}^{h_{\epsilon_j,i+1}(\tau)}\eta(u(x,\tau)|v_{i+1}(x,\tau))\,dx-\sum_{i=0}^{K+1} \int\limits_{h_{\epsilon_j,i}(t)}^{h_{\epsilon_j,i+1}(t)}\eta(u(x,t)|v_{i+1}(x,t))\,dx\\ 
& \leq \sum_{i=0}^{K+1} \int\limits_{t}^{\tau}\bigg[ q(u(h_{\epsilon_j,i}(r)+,r);v_{i+1}(h_{\epsilon_j,i}(r),r))-q(u(h_{\epsilon_j,i+1}(r)-,r);v_{i+1}(h_{\epsilon_j,i+1}(r),r))\\
&\hspace{1in}+\dot{h}_{\epsilon_j,i+1}(r)\eta(u(h_{\epsilon_j,i+1}(r)-,r)|v_{i+1}(h_{\epsilon_j,i+1}(r),r))\\
&\hspace{1in}-\dot{h}_{\epsilon_j,i}(r)\eta(u(h_{\epsilon_j,i}(r)+,r)|v_{i+1}(h_{\epsilon_j,i}(r),r))\bigg]\,dr
\shortintertext{Then, we collect the terms corresponding to $h_{\epsilon_j,i}$ into one sum, and the terms corresponding to $h_{\epsilon_j,i+1}$ into another sum,}
&=\sum_{i=0}^{K+1} \int\limits_{t}^{\tau}\bigg[ q(u(h_{\epsilon_j,i}(r)+,r);v_{i+1}(h_{\epsilon_j,i}(r),r))-\dot{h}_{\epsilon_j,i}(r)\eta(u(h_{\epsilon_j,i}(r)+,r)|v_{i+1}(h_{\epsilon_j,i}(r),r))\bigg]\,dr\\
&\hspace{1in}+\sum_{i=0}^{K+1} \int\limits_{t}^{\tau}\bigg[-q(u(h_{\epsilon_j,i+1}(r)-,r);v_{i+1}(h_{\epsilon_j,i+1}(r),r))\\
&\hspace{1in}+\dot{h}_{\epsilon_j,i+1}(r)\eta(u(h_{\epsilon_j,i+1}(r)-,r)|v_{i+1}(h_{\epsilon_j,i+1}(r),r))\bigg]\,dr
\shortintertext{Next, we peel off the $i=0$ term from the first sum, and the $i=K+1$ term from the second sum,}
&=\int\limits_{t}^{\tau}\bigg[ q(u(h_{\epsilon_j,0}(r)+,r);v_1(h_{\epsilon_j,0}(r),r))-\dot{h}_{\epsilon_j,0}(r)\eta(u(h_{\epsilon_j,0}(r)+,r)|v_1(h_{\epsilon_j,0}(r),r))\bigg]\,dr\\
&+\sum_{i=1}^{K+1} \int\limits_{t}^{\tau}\bigg[ q(u(h_{\epsilon_j,i}(r)+,r);v_{i+1}(h_{\epsilon_j,i}(r),r))-\dot{h}_{\epsilon_j,i}(r)\eta(u(h_{\epsilon_j,i}(r)+,r)|v_{i+1}(h_{\epsilon_j,i}(r),r))\bigg]\,dr\\
&\hspace{1in}+\sum_{i=0}^{K} \int\limits_{t}^{\tau}\bigg[-q(u(h_{\epsilon_j,i+1}(r)-,r);v_{i+1}(h_{\epsilon_j,i+1}(r),r))\\
&\hspace{1in}+\dot{h}_{\epsilon_j,i+1}(r)\eta(u(h_{\epsilon_j,i+1}(r)-,r)|v_{i+1}(h_{\epsilon_j,i+1}(r),r))\bigg]\,dr\\
&\hspace{1in}+\int\limits_{t}^{\tau}\bigg[-q(u(h_{\epsilon_j,K+2}(r)-,r);v_{K+2}(h_{\epsilon_j,K+2}(r),r))\\
&\hspace{1in}+\dot{h}_{\epsilon_j,K+2}(r)\eta(u(h_{\epsilon_j,K+2}(r)-,r)|v_{K+2}(h_{\epsilon_j,K+2}(r),r))\bigg]\,dr
\shortintertext{We then reindex the second sum $\sum_{i=0}^{K}\int\limits_{t}^{\tau}[\cdots]\,dr$ to start at $i=1$, and combine it with the first sum $\sum_{i=1}^{K+1}\int\limits_{t}^{\tau}[\cdots]\,dr$,}
&=\int\limits_{t}^{\tau}\bigg[ q(u(h_{\epsilon_j,0}(r)+,r);v_1(h_{\epsilon_j,0}(r),r))-\dot{h}_{\epsilon_j,0}(r)\eta(u(h_{\epsilon_j,0}(r)+,r)|v_1(h_{\epsilon_j,0}(r),r))\bigg]\,dr\\
&+\sum_{i=1}^{K+1} \int\limits_{t}^{\tau}\bigg[ q(u(h_{\epsilon_j,i}(r)+,r);v_{i+1}(h_{\epsilon_j,i}(r),r))-q(u(h_{\epsilon_j,i}(r)-,r);v_{i}(h_{\epsilon_j,i}(r),r))\\
&-\dot{h}_{\epsilon_j,i}(r)(\eta(u(h_{\epsilon_j,i}(r)+,r)|v_{i+1}(h_{\epsilon_j,i}(r),r))-\eta(u(h_{\epsilon_j,i}(r)-,r)|v_{i}(h_{\epsilon_j,i}(r),r)))\bigg]\,dr\\
&\hspace{1in}+\int\limits_{t}^{\tau}\bigg[-q(u(h_{\epsilon_j,K+2}(r)-,r);v_{K+2}(h_{\epsilon_j,K+2}(r),r))\\
&\hspace{1in}+\dot{h}_{\epsilon_j,K+2}(r)\eta(u(h_{\epsilon_j,K+2}(r)-,r)|v_{K+2}(h_{\epsilon_j,K+2}(r),r))\bigg]\,dr \\
&\leq \frac{\epsilon_j}{T(K+1)}(\tau-t)(K+1) < \epsilon_j
\end{align*}
by \eqref{dissipation_rate3}, the definition of $s$, and noting that $\dot{h}_{\epsilon_j,0}=s$ and $\dot{h}_{\epsilon_j,K+2}=-s$.
Thus, 
\begin{equation}
\begin{aligned}
\label{relative_entropy_decreasing_claim2}
\sum_{i=0}^{K+1} \int\limits_{h_{\epsilon_j,i}(\tau)}^{h_{\epsilon_j,i+1}(\tau)}\eta(u(x,\tau)|v_{i+1}(x,\tau))\,dx
<\\
\sum_{i=0}^{K+1} \int\limits_{h_{\epsilon_j,i}(t)}^{h_{\epsilon_j,i+1}(t)}\eta(u(x,t)|v_{i+1}(x,t))\,dx+ \epsilon_j.
\end{aligned}
\end{equation}

Then, let $j\to\infty$ in \eqref{relative_entropy_decreasing_claim2} and use the dominated convergence theorem to get 
\begin{align}
\label{relative_entropy_decreasing_claim55}
\sum_{i=0}^{K+1} \int\limits_{h_{i}(\tau)}^{h_{i+1}(\tau)}\eta(u(x,\tau)|v_{i+1}(x,\tau))\,dx
&\leq 
\sum_{i=0}^{K+1} \int\limits_{h_{i}(t)}^{h_{i+1}(t)}\eta(u(x,t)|v_{i+1}(x,t))\,dx
\end{align}
for almost every $t$ and $\tau$ verifying $t^*\leq t < \tau<t^{**}$.

From \eqref{relative_entropy_decreasing_claim55}, we get
\begin{align}
\label{relative_entropy_decreasing_claim555}
{\rm ap}\,\lim_{\tau\to {t^{**}}^{-}}\sum_{i=0}^{K+1} \int\limits_{h_{i}(t^{**})}^{h_{i+1}(t^{**})}\eta(u(x,\tau)|v_{i+1}(x,\tau))\,dx\leq\\
{\rm ap}\,\lim_{t\to {t^{*}}^{+}}\sum_{i=0}^{K+1} \int\limits_{h_{i}(t^{*})}^{h_{i+1}(t^{*})}\eta(u(x,t)|v_{i+1}(x,t))\,dx
\end{align}
where we have used

\begin{align}
\label{fixed_time_terminals_ap_lim_a}
{\rm ap}\,\lim_{\tau\to {t^{**}}^{-}}\sum_{i=0}^{K+1} \int\limits_{h_{i}(\tau)}^{h_{i+1}(\tau)}\eta(u(x,\tau)|v_{i+1}(x,\tau))\,dx=\\
{\rm ap}\,\lim_{\tau\to {t^{**}}^{-}}\sum_{i=0}^{K+1} \int\limits_{h_{i}(t^{**})}^{h_{i+1}(t^{**})}\eta(u(x,\tau)|v_{i+1}(x,\tau))\,dx
\end{align}
and
\begin{align}
\label{fixed_time_terminals_ap_lim_b}
{\rm ap}\,\lim_{t\to {t^{*}}^{+}}\sum_{i=0}^{K+1} \int\limits_{h_{i}(t)}^{h_{i+1}(t)}\eta(u(x,t)|v_{i+1}(x,t))\,dx={\rm ap}\,\lim_{t\to {t^{*}}^{+}}\sum_{i=0}^{K+1} \int\limits_{h_{i}(t^{*})}^{h_{i+1}(t^{*})}\eta(u(x,t)|v_{i+1}(x,t))\,dx.\end{align}

The approximate limits exist by \Cref{left_right_ap_limits}. Then \eqref{relative_entropy_decreasing_claim555} and \eqref{left_and_right_limits_order} give the claim \eqref{relative_entropy_decreasing_claim1}.
\end{claimproof}

If $t^{**}=T$, then we have proven \Cref{psi_exists_with_limits} holds for $N=K+1$: define $x_{i,T}\coloneqq h_{i}(T)$ for $0\leq i \leq K+2$. By \eqref{relative_entropy_decreasing_claim1}, the $x_{0,T},x_{1,T},\ldots,x_{K+2,T}$ satisfy the conclusions of \Cref{psi_exists_with_limits}. 

Otherwise, $t^{**}<T$ and we consider the $0\leq i,j\leq K+3$  such that the following holds:
\begin{align}
\begin{cases} \label{good_i_and_j}
h_{i-1}(t^{**})<h_{i}(t^{**})<h_{j}(t^{**})\\
\mbox{for $i<k<j$, $h_{i}(t^{**})=h_{k}(t^{**})$}
\end{cases}
\end{align}
where $h_{-1}(t)\coloneqq-\infty$ and $h_{K+3}(t)\coloneqq+\infty$ for all $t$. Then, let $\{(i_n,j_n)\}_{n\in\{1,\ldots,L\}}$ for $L\in\mathbb{Z}_{+}$ be the set of $i$ and $j$ pairs which satisfy \eqref{good_i_and_j} (label the $i$ and $j$ pairs such that $i_n < i_{n+1}$ for all $n$). Note that each $i$ has only one corresponding $j$. Thus, due to at least two of the $h_{i}(t^{**})$ equalling each other (for $i$ ranging over $0,\ldots,K+2$), $L\leq K+2$. 

By the induction hypothesis with $N=L-2$: there exist real numbers $\tilde{x}_{0,T},\ldots,\tilde{x}_{L-1,T}$ verifying $\tilde{x}_{0,T}=-R\leq\tilde{x}_{1,T}\leq \cdots\leq \tilde{x}_{L-2,T}\leq R=\tilde{x}_{L-1,T}$  and
\begin{equation}
\begin{aligned}
\label{l2_stability_induction_hypothesis}
{\rm ap}\,\lim_{t\to {T}^{+}}\sum_{l=0}^{L-2} \int\limits_{\tilde{x}_{l,T}}^{\tilde{x}_{l+1,T}}\eta(u(x,t)|v_{i_{l+2}}(x,t))\,dx
\leq \\
{\rm ap}\,\lim_{t\to {t^{**}}^{+}}\sum_{l=0}^{L-2} \int\limits_{h_{i_{l+1}}(t^{**})}^{h_{i_{l+2}}(t^{**})}\eta(u(x,t)|v_{i_{l+2}}(x,t))\,dx.
\end{aligned}
\end{equation}

For each $n\in\{1,\ldots,L\}$, define
\begin{align}
x_{i,T}\coloneqq \tilde{x}_{n-1,T} \mbox{ for all $i_n\leq i<j_n$.}
\end{align}

Then by construction $x_{0,T},\ldots,x_{K+2,T}$ satisfy the conclusions of \Cref{psi_exists_with_limits} with $N=K+1$. In particular, \eqref{l2_stability1} follows from \eqref{relative_entropy_decreasing_claim1} and \eqref{l2_stability_induction_hypothesis}.

Thus, by the principle of strong induction we have proven \Cref{psi_exists_with_limits} for all $N\in\mathbb{N}\cup\{0\}$.

\subsection{Proof of \Cref{LIWAS_dense}}
Step functions are dense in $L^2([-M,M])$. Thus, there exists a step function $s \in L^2([-M,M])$ such that $\norm{f-s}_{L^2([-M,M])} < \frac{\epsilon}{2}$ and $\norm{s}_{L^\infty([-M,M])} \leq \norm{f}_{L^\infty([-M,M])}$. We can write $s$ in the form
\begin{align}
s(x)=\bigg(\sum_{i=1}^{n_{+}} \alpha^{+}_i H(x-x^{+}_{i})\bigg)+\bigg(\sum_{i=1}^{n_{-}} \alpha^{-}_i H(x-x^{-}_{i})\bigg)
\end{align}
for some $n_{+},n_{-}\in\mathbb{N}$, $\{\alpha^{+}_i\}_{i=1}^{n_{+}}\subset(0,\infty)$, $\{\alpha^{-}_i\}_{i=1}^{n_{-}}\subset(-\infty,0)$, $\{x^{+}_i\}_{i=1}^{n_{+}}\subset(-M,M)$ and $\{x^{-}_i\}_{i=1}^{n_{-}}\subset(-M,M)$. $H$ is the Heaviside step function
\begin{align}
H(x)\coloneqq
\begin{cases}
0 & \text{if } x<0\\
1 & \text{if } x>0.
\end{cases}
\end{align}

Define
\begin{align}
s^{+}(x)\coloneqq \sum_{i=1}^{n_{+}} \alpha^{+}_i H(x-x^{+}_{i}),\\
s^{-}(x)\coloneqq \sum_{i=1}^{n_{-}} \alpha^{-}_i H(x-x^{-}_{i}).
\end{align}
We can then write $s=s^{+}+s^{-}$.

Consider the standard mollifier $m: \mathbb{R}\to\mathbb{R}$, where $m$ is smooth and compactly supported, $m\geq0$, and $\int m(x)\,dx=1$.

Let $\delta>0$. Define 
\begin{align}
m_\delta(x)\coloneqq \frac{1}{\delta} m\bigg(\frac{x}{\delta}\bigg).
\end{align}

Define
\begin{align}
f_\delta (x)\coloneqq 
\begin{cases}
\displaystyle\lim_{y\to -M^{+}}[(s^{+}* m_\delta)(y)+s^{-}(y)] & \text{if } x\leq-M\\
(s^{+}* m_\delta)(x)+s^{-}(x) & \text{if } -M<x<M\\
\displaystyle\lim_{y\to M^{-}}[(s^{+}* m_\delta)(y)+s^{-}(y)] & \text{if } x\geq M.
\end{cases}
\end{align}

Note that $f_\delta$ is of the form \eqref{LIWAS}, $\norm{f_\delta}_{L^\infty(\mathbb{R})} \leq \norm{f}_{L^\infty([-M,M])}$ when $\delta < \inf_{i,j}|x^{+}_i-x^{-}_j|$, and all of the discontinuities in $f_\delta$ are in the interval $(-M,M)$ because $\{x^{-}_i\}_{i=1}^{n_{-}}\subset(-M,M)$.

From the Minkowski inequality, 
\begin{equation}
\begin{aligned}
\label{triangle_ineq_LIWAS}
\norm{f-f_\delta}_{L^2([-M,M])}&\leq \norm{f-s}_{L^2([-M,M])} + \norm{s-f_\delta}_{L^2([-M,M])}\\
&< \frac{\epsilon}{2} + \norm{(s^{+}* m_\delta)-s^{+}}_{L^2([-M,M])}.
\end{aligned}
\end{equation}

Choose $\delta$ even smaller such that $\norm{(s^{+}* m_\delta)-s^{+}}_{L^2([-M,M])}< \frac{\epsilon}{2}$. This completes the proof.

\section{Main proposition implies main theorem}\label{main_prop_implies_main_theorem}
Let $u\in L^{\infty}(\mathbb{R}\times[0,\infty))$ be a weak solution to \eqref{conservation_law}, with initial data $u^0$. Let $\eta\in C^2(\mathbb{R})$ be a strictly convex entropy. Assume that $u$ is entropic for the entropy $\eta$ and verifies the strong trace property (\Cref{strong_trace_condition}).

We will show that $u$ satisfies \eqref{condition_e}.

First, by strict convexity of $\eta$ we can choose a constant $c^*>0$ such that
\begin{align}\label{control_rel_entropy_below1}
c^{*}(a-b)^2\leq \eta(a|b)
\end{align}
for all $a,b\in[-B,B]$, where $B$ is defined in \eqref{define_B}.

Let $R,T>0$. From \Cref{main_lemma}, for all $\epsilon>0$, there exists $\psi_\epsilon: [-R,R]\to\mathbb{R}$ such that
\begin{align}\label{relative_entropy_stable_epsilon}
\int\limits_{\abs{x}\leq R}\eta(u(x,T)|\psi_\epsilon(x))\,dx \leq c^* \epsilon^2
\end{align}
and
 \begin{align}\label{psi_satisfies_condition_e_epsilon}
\psi_\epsilon(x+z)-\psi_\epsilon(x)\leq \frac{c}{T}z
\end{align}
for $x\in[-R,R]$ and $z>0$ with $x+z\in[-R,R]$. Here $c=1/\inf{A''}$. We have also 
\begin{align}\label{psi_bounded_u_epsilon}
\norm{\psi_\epsilon}_{L^\infty([-R,R])}\leq \norm{u}_{L^\infty(\mathbb{R}\times[0,\infty))}.
\end{align}

By \eqref{psi_bounded_u_epsilon}, $\norm{\psi_\epsilon}_{L^\infty([-R,R])}\leq B$. Likewise, we have $\norm{u(\cdot,T)}_{L^\infty(\mathbb{R})}\leq B$. Thus, from \eqref{control_rel_entropy_below1}, we get
\begin{align}\label{l_2_control_psi_epsilon_u31231}
c^*(\psi_\epsilon(x)-u(x,T))^2\leq \eta(u(x,T)|\psi_\epsilon(x))
\end{align}
for almost every $x\in\ [-R,R]$.

Then from \eqref{relative_entropy_stable_epsilon} and \eqref{l_2_control_psi_epsilon_u31231} we have
\begin{align}\label{l_2_control_psi_epsilon_u}
\norm{\psi_\epsilon(\cdot)-u(\cdot,T)}_{L^2([-R,R])}\leq \epsilon.
\end{align}

Thus, there exists a sequence $\{\epsilon_j\}_{j=1}^{\infty}$ with $\epsilon_j\to0^{+}$ such that
\begin{align}
\psi_{\epsilon_j}(x)\to u(x,T) \mbox{ as $j\to\infty$ for almost every x.}
\end{align}

Additionally, from \eqref{psi_satisfies_condition_e_epsilon}:
\begin{align}\label{psi_j_satisfies_condition_e_epsilon}
\psi_{\epsilon_j} (x+z)-\psi_{\epsilon_j}(x)\leq \frac{c}{T}z
\end{align}
for all $j\in\mathbb{N}$, $x\in[-R,R]$, and $z>0$ with $x+z\in[-R,R]$. 

We let $j\to\infty$ in \eqref{psi_j_satisfies_condition_e_epsilon} to get:
\begin{align}\label{u_satisfies_condition_from_psi_j}
u(x+z,T)-u(x,T)\leq \frac{c}{T}z
\end{align}
for almost every $x\in[-R,R]$, and almost every $z>0$ with $x+z\in[-R,R]$. 

Because $R,T>0$ are arbitrary, \eqref{u_satisfies_condition_from_psi_j} implies that $u$ satisfies \eqref{condition_e}. This concludes the proof of \Cref{main_theorem}.

\section{Appendix}
\subsection{Proof of \Cref{d_sm_rh_lemma}} \label{d_sm_rh_lemma_proof}

Let $u_L,u_R,u_{-},u_{+}\in\mathbb{R}$. We then borrow the following notation from \cite{serre_vasseur}: If $F$ is a function of $u$, then we define
\begin{align}
F_{R}\coloneqq F(u_{R}), \hspace{.1in}F_{L}\coloneqq F(u_{L}), \hspace{.1in} [F] \coloneqq F_{R}-F_{L}, \hspace{.1in} F_{\pm}\coloneqq F(u_\pm).
\end{align}

\emph{Proof of \eqref{d_rh}}

This proof is from \cite[p.~9-10]{serre_vasseur}. In \cite{serre_vasseur}, for the general systems case, the authors develop a condition which they label with the equation number 7 \cite[p.~4]{serre_vasseur}. In particular, they claim to use this condition to show  \eqref{d_rh} for the scalar case. In fact, the condition is not necessary in the scalar case and their proof goes through unchanged without the condition.

Denote
\begin{align}
D\coloneqq q(u_{+};u_R)-q(u_{-};u_L)-\sigma(u_{-},u_{+})(\eta(u_{+}|u_R)-\eta(u_{-}|u_L)).
\end{align}
Further, let $\sigma$ denote $\sigma(u_{-},u_{+})$.

From Rankine-Hugoniot (as noted in \cite[p.~5]{serre_vasseur}), 
\begin{align}
D=[\eta'(u)A(u)-q(u)]-\sigma[\eta'(u)u-\eta]+q_{+}-q_{-}-\sigma(\eta_{+}-\eta_{-})-[\eta'](A(u)-\sigma u)_{\pm},
\end{align}
where $(A(u)-\sigma u)_{\pm}$ denotes that 
\begin{align}
A(u_{+})-\sigma u_{+}=A(u_{-})-\sigma u_{-} \label{equality_RH1}
\end{align}
because of the Rankine-Hugoniot relation \eqref{rh_velocity}.

From the fundamental theorem of calculus and integration by parts,
\begin{align}
D=\int\limits_{u_{+}}^{u_{-}} \eta''(u)(A(u)-\sigma u)\,du -\int\limits_{u_{+}}^{u_{-}} \eta''(u)\,du (A(u)-\sigma u)_{\pm}\\
-\int\limits_{u_{R}}^{u_{L}} \eta''(u)(A(u)-\sigma u)\,du+\int\limits_{u_{R}}^{u_{L}} \eta''(u)\,du (A(u)-\sigma u)_{\pm}.
\end{align}
We can then write
\begin{align}
D=\epsilon(I)B(I)+\epsilon(J)B(J)
\end{align}
where $I, J$ are disjoint intervals such that
\begin{align}
I\cup J =((u_{+},u_{-})\cup(u_R,u_L))\setminus ((u_{+},u_{-})\cap(u_R,u_L)).
\end{align}
We define the sign $\epsilon(I)$ to be $+1$ if $I\subset (u_{+},u_{-})$ and $-1$ otherwise. Finally,
\begin{align}
B(I)\coloneqq \int\limits_{I}\eta''(u)(A(u)-\sigma u)\,du - (A(u)-\sigma u)_{\pm} \int\limits_{I}\eta''(u)\,du.
\end{align}
The function $u\mapsto A(u)-\sigma u$ is strictly convex and \eqref{equality_RH1} holds. Thus, $B(I)<0$ if $I\subset (u_{+},u_{-})$ and positive otherwise. Thus, for all intervals $\epsilon(I)B(I)<0$. Thus $D<0$.

\emph{Proof that \eqref{denom_zero} implies \eqref{d_sm}}

This proof is from \cite[p.~9]{serre_vasseur}.

Remark that the equality 
\begin{align}
\eta(u|u_L)=\eta(u|u_R)
\end{align}
is equivalent to \cite[p.~4]{serre_vasseur}
\begin{align}
[\eta']u=[\eta'(u)u-\eta(u)].
\end{align}
Remark also (as noted in \cite[p.~4]{serre_vasseur})
\begin{align}
q(u;u_R)-q(u;u_L)=[\eta'A-q]-[\eta']A(u).
\end{align}

Then, \eqref{d_sm} is equivalent to (as noted in \cite[p.~9]{serre_vasseur})
\begin{align}
\frac{\int\limits_{u_R}^{u_L} \eta''(u)A(u)\,du}{\int\limits_{u_R}^{u_L}\eta''(u)\,du} > A\Bigg(\frac{\int\limits_{u_R}^{u_L} \eta''(u)u\,du}{\int\limits_{u_R}^{u_L}\eta''(u)\,du} \Bigg). \label{equivalent1}
\end{align}
Finally, \eqref{equivalent1} is true by Jensen's inequality because $\eta''>0$ so $\eta''(u)\,du$ is a measure, and $A$ is strictly convex.

\subsection{Proof of \Cref{Filippov20}}\label{proof_of_Filippov}

The following proof of \eqref{Filippov1}, \eqref{Filippov2} and \eqref{Filippov3} is based on the proof of Proposition 1 in \cite{Leger2011} and the proof of Lemma 2.2 in \cite{serre_vasseur}. We do not prove \eqref{Filippov4} and \eqref{Filippov5}: these properties are in Lemma 6 in \cite{Leger2011}, and their proofs are in the appendix in \cite{Leger2011}.

Define
\begin{align}
v_n(x,t)\coloneqq \int\limits_{0}^{1} V_\epsilon(u(x+\frac{y}{n},t),\bar{u}_1(x+\frac{y}{n},t),\bar{u}_2(x+\frac{y}{n},t))\,dy.
\end{align}

Let $h_{\epsilon,n}$ be the solution to the ODE:
\begin{align}
\begin{cases}\label{ode_for_h_epsilon_n}
\dot{h}_{\epsilon,n}(t)=v_n(h_{\epsilon,n}(t),t), \text{ for } t>0\\
h_{\epsilon,n}(t^*)=x_0.
\end{cases}
\end{align}
Due to $V_\epsilon$ being continuous, $v_n$ is bounded uniformly in $n$ ($\hspace{.03in}\norm{v_n}_{L^{\infty}}\leq \norm{V_\epsilon}_{L^{\infty}}$), Lipschitz continuous in $x$, and measurable in $t$. Thus \eqref{ode_for_h_epsilon_n} has a unique solution in the sense of Carath\'eodory.  

Because $v_n$ is bounded uniformly in $n$, the $h_{\epsilon,n}$ are Lipschitz continuous in time, with the Lipschitz constant uniform in $n$. Thus, by Arzel\`a--Ascoli the $h_{\epsilon,n}$ converge in $C^0(0,T)$ for any fixed $T>0$ to a Lipschitz continuous function $h_{\epsilon}$ (passing to a subsequence if necessary). Note that $\dot{h}_{\epsilon,n}$ converges in $L^{\infty}$ weak* to $\dot{h}_{\epsilon}$.

We define
\begin{align}
V_{\text{max}}(t)\coloneqq \max\{V_\epsilon(u_{-},\bar{u}_1,\bar{u}_2),V_\epsilon(u_{+},\bar{u}_1,\bar{u}_2)\},\\
V_{\text{min}}(t)\coloneqq \min\{V_\epsilon(u_{-},\bar{u}_1,\bar{u}_2),V_\epsilon(u_{+},\bar{u}_1,\bar{u}_2)\},
\end{align}
where $u_{\pm}=u(h_\epsilon(t)\pm,t)$ and $\bar{u}_i=\bar{u}_i(h_\epsilon(t),t)$.

To show \eqref{Filippov3}, we will first prove that for almost every $t>0$
\begin{align}
\lim_{n\to\infty}[\dot{h}_{\epsilon,n}(t)-V_{\text{max}}(t)]_{+}=0, \label{one_side_filippov}\\
\lim_{n\to\infty}[V_{\text{min}}(t)-\dot{h}_{\epsilon,n}(t)]_{+}=0,\label{other_side_filippov}
\end{align}
where $[\hspace{.025in}\cdot\hspace{.025in}]_{+}\coloneqq \max(0,\cdot)$.

The proofs of \eqref{one_side_filippov} and \eqref{other_side_filippov} are similar. We will only show the first one.

\begin{align}
[\dot{h}_{\epsilon,n}(t)-&V_{\text{max}}(t)]_{+}\\
&=\Bigg[\int\limits_{0}^{1} V_\epsilon(u(h_{\epsilon,n}(t)+\frac{y}{n},t),\bar{u}_1(h_{\epsilon,n}(t)+\frac{y}{n},t),\bar{u}_2(h_{\epsilon,n}(t)+\frac{y}{n},t))\,dy-V_{\text{max}}(t)\Bigg]_{+}\\
&=\Bigg[\int\limits_{0}^{1} V_\epsilon(u(h_{\epsilon,n}(t)+\frac{y}{n},t),\bar{u}_1(h_{\epsilon,n}(t)+\frac{y}{n},t),\bar{u}_2(h_{\epsilon,n}(t)+\frac{y}{n},t))-V_{\text{max}}(t)\,dy\Bigg]_{+}\\
&\leq\int\limits_{0}^{1} [V_\epsilon(u(h_{\epsilon,n}(t)+\frac{y}{n},t),\bar{u}_1(h_{\epsilon,n}(t)+\frac{y}{n},t),\bar{u}_2(h_{\epsilon,n}(t)+\frac{y}{n},t))-V_{\text{max}}(t)]_{+}\,dy\\
&\leq \esssup_{y\in(0,\frac{1}{n})}[V_\epsilon(u(h_{\epsilon,n}(t)+y,t),\bar{u}_1(h_{\epsilon,n}(t)+y,t),\bar{u}_2(h_{\epsilon,n}(t)+y,t))-V_{\text{max}}(t)]_{+}\\
&\leq \esssup_{y\in(-\epsilon_n,\epsilon_n)}[V_\epsilon(u(h_{\epsilon}(t)+y,t),\bar{u}_1(h_{\epsilon}(t)+y,t),\bar{u}_2(h_{\epsilon}(t)+y,t))-V_{\text{max}}(t)]_{+},\label{final_line_for_V_max_bound_proof_of_Filippov}
\end{align}
where $\epsilon_n\coloneqq \abs{h_{\epsilon,n}(t)-h_{\epsilon}(t)}+\frac{1}{n}$. Note $\epsilon_n\to0^{+}$.

Fix a $t\geq0$ such that we have a strong trace according to \Cref{strong_trace_condition}. Then, by the continuity of $V_\epsilon$,
\begin{align}
\lim_{n\to\infty}\esssup_{y\in(0,\frac{1}{n})}[V_\epsilon(u(h_{\epsilon}(t)\pm y,t),\bar{u}_1(h_{\epsilon}(t)\pm y,t),\bar{u}_2(h_{\epsilon}(t)\pm y,t))-V_\epsilon(u_{\pm},\bar{u}_1,\bar{u}_2)]_{+}=0,
\end{align}
where $u_{\pm}=u(h_\epsilon(t)\pm,t)$ and $\bar{u}_i=\bar{u}_i(h_\epsilon(t),t)$.

This implies 
\begin{align}
\lim_{n\to\infty}\esssup_{y\in(0,\frac{1}{n})}[V_\epsilon(u(h_{\epsilon}(t)\pm y,t),\bar{u}_1(h_{\epsilon}(t)\pm y,t),\bar{u}_2(h_{\epsilon}(t)\pm y,t))-V_{\text{max}}(t)]_{+}=0.\label{thing_we_need_goes_to_zero_proof_of_Filippov}
\end{align}

We can control \eqref{final_line_for_V_max_bound_proof_of_Filippov} from above by
\begin{equation}
\begin{aligned}\label{final_line_for_V_max_bound_proof_of_Filippov_next_line}
\esssup_{y\in(-\epsilon_n,0)}[V_\epsilon(u(h_{\epsilon}(t)+y,t),\bar{u}_1(h_{\epsilon}(t)+y,t),\bar{u}_2(h_{\epsilon}(t)+y,t))-V_{\text{max}}(t)]_{+}+\\
\esssup_{y\in(0,\epsilon_n)}[V_\epsilon(u(h_{\epsilon}(t)+y,t),\bar{u}_1(h_{\epsilon}(t)+y,t),\bar{u}_2(h_{\epsilon}(t)+y,t))-V_{\text{max}}(t)]_{+}.
\end{aligned}
\end{equation}

Then, by \eqref{thing_we_need_goes_to_zero_proof_of_Filippov}, the quantity \eqref{final_line_for_V_max_bound_proof_of_Filippov_next_line} goes to $0$ as $n\to\infty$. This proves \eqref{one_side_filippov}.

Recall that $\dot{h}_{\epsilon,n}$ converges in $L^{\infty}$ weak* to $\dot{h}_{\epsilon}$. Thus, because the function $[\hspace{.025in}\cdot\hspace{.025in}]_{+}$ is convex, 
\begin{align}
 \int\limits_{0}^{T} [\dot{h}_{\epsilon}(t)-V_{\text{max}}(t)]_{+} \,dt \leq \liminf_{n\to\infty} \int\limits_{0}^{T}[\dot{h}_{\epsilon,n}(t)-V_{\text{max}}(t)]_{+}\,dt.
\end{align}

By the dominated convergence theorem and \eqref{one_side_filippov},
\begin{align}
\liminf_{n\to\infty} \int\limits_{0}^{T}[\dot{h}_{\epsilon,n}(t)-V_{\text{max}}(t)]_{+}\,dt=0.
\end{align}

This proves 
\begin{align}
 \int\limits_{0}^{T} [\dot{h}_{\epsilon}(t)-V_{\text{max}}(t)]_{+} \,dt=0.
\end{align}
 
 A similar argument gives
 \begin{align}
 \int\limits_{0}^{T} [V_{\text{min}}(t)-\dot{h}_{\epsilon}(t)]_{+} \,dt=0.
 \end{align}
 
 This proves \eqref{Filippov3}.
 
\bibliographystyle{plain}
\bibliography{references}
\end{document}